%% file: main-arxiv.tex
\documentclass[reqno,letterpaper]{amsart}

\RequirePackage[OT1]{fontenc}
\usepackage[tchypersetup]{tchdr}
\usepackage{natbib}

\numberwithin{equation}{section}

\newcommand{\restrict}[2]{ {\left. { #1 } \right|_{#2}} }
\newcommand{\order}[1]{ {\left[ { #1 } \right]} }
\newcommand{\memb}[2]{ { #2^{(#1)} } }
\newcommand{\mult}[1]{ {\overline{ { #1 } }} }
\newcommand{\traits}{\mathbb{T}}
\newcommand{\squash}{\mathrm{filter}}
\newcommand{\labels}{\mathbb{Y}}

\newcommand{\memseqs}{\mathbb{K}}

\begin{document}
\sloppy
\allowdisplaybreaks

\title{Exchangeable Trait Allocations}
\author[T.~Campbell]{Trevor Campbell}
\address{Computer Science and Artificial Intelligence Laboratory (CSAIL)\\ Massachusetts Institute of Technology}
\urladdr{http://www.trevorcampbell.me/}
\email{tdjc@mit.edu}

\author[D.~Cai]{Diana Cai}
\address{Department of Computer Science\\Princeton University}
\urladdr{http://www.dianacai.com/}
\email{dcai@cs.princeton.edu}

\author[T.~Broderick]{Tamara Broderick}
\address{Computer Science and Artificial Intelligence Laboratory (CSAIL) \\ Massachusetts Institute of Technology}
\urladdr{http://www.tamarabroderick.com}
\email{tbroderick@csail.mit.edu}


\input{abstract}
\maketitle
\input{intro}

\input{conventions}

\input{traitalloc}

\input{orderofappearance}
\input{paintbox}

\input{freq}

\input{applications}

\input{conc}

\appendix
\input{proofs}

\newpage
\small
\bibliographystyle{ba}
\bibliography{main}

\end{document}

%% file: abstract.tex
\begin{abstract}

Trait allocations are a class of combinatorial structures in which data
may belong to multiple groups and may have different levels of belonging in each group.
Often the data are also exchangeable, i.e., their joint distribution 
is invariant to reordering. In clustering---a special case 
of trait allocation---exchangeability implies the existence of both a de Finetti representation
and an exchangeable partition probability function (EPPF),  
distributional representations useful for computational and theoretical purposes. 
In this work, we develop the analogous de Finetti representation and exchangeable trait probability function (ETPF)
for trait allocations, along with a characterization of all trait allocations with an ETPF. Unlike previous 
feature allocation characterizations, our proofs fully capture single-occurrence ``dust'' groups.
We further introduce a novel constrained version of the ETPF that we use to establish an intuitive connection
between the probability functions for clustering, feature allocations, and trait allocations.
As an application of our general theory, we characterize the distribution of all edge-exchangeable graphs, 
a class of recently-developed models that captures realistic sparse graph sequences.


\end{abstract}

%% file: intro.tex
\section{Introduction}

Representation theorems for exchangeable random variables are a ubiquitous and powerful tool in Bayesian modeling
and inference. In many data analysis problems, we impose an order, or indexing, on our data points.
This indexing can arise naturally---if we are truly observing data in a sequence---or can be artificially 
created to allow their storage in a database. In this context, exchangeability expresses the assumption that this order
is arbitrary and should not affect our analysis.
For instance, we often assume a sequence of data points 
is an \emph{infinite exchangeable sequence},
i.e., that the distribution
of any finite subsequence is invariant to reordering.
Though this assumption may seem weak, de Finetti's theorem
\citep{deFinetti31,Hewitt55}
tells us that in this case, we can assume that a latent parameter exists, that our data
are independent and identically distributed (\iid) conditional on this parameter,
and that the parameter itself has a distribution. Thus, de Finetti's theorem may be
seen as a justification for a Bayesian model and prior---and, in fact, for the
infinite-dimensional priors provided by Bayesian nonparametrics \citep{Jordan10}.

De Finetti-style representation theorems have provided many other useful insights for modeling and
inference within Bayesian analysis. For example, consider \emph{clustering} problems, where the inferential goal is to assign data
points to mutually exclusive and exhaustive groups.
It is typical to assume that the distribution of the clustering---i.e., the assignment of
data points to clusters---is invariant
to the ordering of the data points. In this case, two different representation theorems have proved particularly useful in practice. First, \citet{Kingman78} showed that exchangeability in clustering implies
the existence of a latent set of probabilities (known as the ``Kingman paintbox'') 
from which cluster assignments are chosen \iid.
It is straightforward to show from the Kingman paintbox representation that exchangeable clustering models
enforce linear growth in cluster size as a function of the size of the total data.
By contrast, many real-world clustering problems, such as disambiguating census data or clustering academic papers by originating lab,
exhibit sublinear growth in cluster size \citep[e.g.,][]{Wallach10,Broderick14,Miller15}.
Thus, the Kingman paintbox representation allows us to see that exchangeable clustering models
are misspecified for these examples. 
Similarly, \citet{Pitman95} showed that clustering exchangeability
is equivalent to the existence of an \emph{exchangeable partition probability function} (EPPF).
The EPPF and similar developments have led to algorithms that allow practical inference
specifically with the Dirichlet process mixture \citep{Escobar94,Escobar95} and more generally in other clustering models
\citep{Pitman97,Ishwaran01,Ishwaran03,Lee13}.

In this work, we develop and characterize a generalization of clustering models that we call \emph{trait allocation models}.
Trait allocations apply when data may belong to more than one group (a \emph{trait}), and may exhibit nonnegative integer levels of belonging in each group.
For example, a document might exhibit multiple words in a number of topics, a participant in a social network
might send multiple messages to each of her friend groups, or a DNA sequence might exhibit different
numbers of genes from different ancestral populations.
Trait allocations generalize both clustering, where data must belong to exactly one group, and feature 
allocations \citep{Griffiths05,Broderick13}, where data exhibit binary membership in multiple groups. 
Authors have recently proposed a number of models for trait allocations \citep[e.g.,][]{Titsias08,Zhou12,Zhou14,James14,Broderick15,Roychowdhury15}.
But as of yet, there is no characterization either of the class of exchangeable
trait allocation models or of classes of exchangeable trait allocation models
that are particularly amenable to inference. The consequences of the exchangeability assumption in 
this setting have not been explored. In this work, we provide characterizations of
both the full class of exchangeable trait allocations and those with EPPF-like probability distributions.
This work not only unifies and generalizes past research on partitions and feature allocations, but
provides a natural avenue for the study of other practical exchangeable combinatorial structures.

We begin by formally defining trait allocations, random sequences thereof, and exchangeability in \cref{sec:traitalloc}.
In \cref{sec:orderofappearance}, we introduce ordered trait allocations via the lexicographic ordering. We use these constructions to establish a de Finetti representation for exchangeable trait allocations in \cref{sec:paintbox}
that is analogous to the Kingman paintbox representation for clustering.
Our new representation handles \emph{dust}, the case where some traits may
appear for just a single data point.  This work therefore also extends previous
work on the special case of exchangeable feature allocations
to the fully general case, whereas
previously it was restricted to the dustless case \citep{Broderick13}.
In \cref{sec:freq}, we develop an EPPF-like function to describe distributions over exchangeable trait allocations
and characterize the class of trait allocations to which it applies.
We call these \emph{exchangeable trait probability functions} (ETPFs).
Just as in the partition and feature allocation cases, the class of random trait allocations with
probability functions represents
a class of trait allocations that are particularly amenable to approximate posterior inference in practice---and
therefore of particularly pressing interest to characterize.
In \cref{sec:freq}, we introduce new concepts we call \emph{constrained ETPFs}, which are the combinatorial
analogue of earlier work on restricted nonparametric processes \citep{Williamson13,DoshiVelez17}.
In \cref{sec:freq,sec:applications}, we show how constrained ETPFs
capture earlier probability functions for numerous exchangeable models within a single framework.
In \cref{sec:applications}, we apply both our de Finetti representation and constrained ETPF
to characterize \emph{edge-exchangeable} graphs, a recently
developed form of exchangeability for graph models that allows sparse projective sequences of graphs \citep{Broderick15a,Crane15,Cai16,Crane16,Williamson16}.
A similar representation generalizing partitions and edge-exchangeable (hyper)graphs
has been studied in concurrent work \citep{Crane16b} on \emph{relational exchangeability}, first introduced by \citet{Ackerman15,Crane15a}---but
here we additionally explore the existence of a trait frequency model, the existence of a constrained trait frequency model
and its connection to clustering and feature allocations, and the various connections between frequency models and probability functions.

%% file: conventions.tex
\subsection{Notation and conventions}\label{sec:conventions}
Definitions are denoted by the symbol $\defined$.
The natural numbers are denoted $\nats \defined \{1, 2, \dots\}$ and the nonnegative reals  $\reals_+ \defined [0, \infty)$.
We let $[N] \defined \{1, 2, \dots, N\}$ for any $N\in\nats$.
Sequences are denoted with parentheses, with indices suppressed
only if they are clear from context. For example, $(x_k)$ is the sequence $(x_k)_{k\in\nats}$
and $(x_{kj})$ is the sequence $(x_{kj})_{k,j\in\nats}$, while $(x_{kj})_{j=1}^\infty$
is the sequence $x_{k1}, x_{k2}, \dots$ with $k$ fixed.
The notation $A\subset B$ means $A$ is a (not necessarily proper) subset of $B$.
The indicator function is denoted $\ind\left(\dots\right)$; for example, $\ind\left(x\in A\right)$ is
1 if $x\in A$, and 0 otherwise.
For any multiset $x$ of elements in a set $\mcX$, we denote $x(y)$
to be the multiplicity of $y$ in $x$ for each $y\in\mcX$.
Two multisets $x, x'$ of $\mcX$ are said to be equal, denoted $x=x'$, if the multiplicity of all elements $y\in\mcX$
are equal in both $x$ and $x'$, i.e.~$\forall y\in\mcX, \, x(y)=x'(y)$.
For any finite or infinite sequence, we use subscript $k$ to denote the $k^\text{th}$ element in the sequence.
For sequences of (multi)sets, if $k$ is beyond the end of the sequence, the subscript $k$ operation returns the empty set.
Equality in distribution and almost surely are denoted $\eqd/\eqas$, and convergence almost surely/in probability/in distribution is denoted $\convas/\convp/\convd$.
We often use cycle notation for permutations~(see \citet[p.~29]{Dummit04}): for example, $\pi = (12)(34)$ is the permutation $\pi$ with $\pi(1)=2$, $\pi(2)=1$, $\pi(3)=4$, $\pi(4)=3$,
and $\pi(k) = k$ for $k>4$.
We use the notation $X \dist (\theta_{j})_{j=0}^\infty$ to denote sampling $X$ from
the categorical distribution on $\{0\}\cup\nats$ with probabilities $\Pr\left(X = j\right) = \theta_{j}$ for $j\in\{0\}\cup\nats$.
The symbol $\bigtimes_N S_N$ for a sequence of sets $(S_N)$ denotes
$S_1\times S_2\times \dots$, their infinite product space.

%% file: traitalloc.tex
\section{Trait allocations}\label{sec:traitalloc}
We begin by formalizing the concepts of a \emph{trait} and \emph{trait allocation}.
We assume that our sequence of data points is indexed by $\nats$.
As a running example for intuition, consider the case where each data point is a document,
and each trait is a topic. Each document may have multiple words that belong to each topic.
The degree of membership of the document in the topic is the number of words in that topic.
We wish to capture the assignment of data points to the traits they express
but in a way that does not depend on the type of data at hand. Therefore, we focus on the
indices to the data points.
This leads to the definition of traits as multisets of the data indices, i.e., the natural numbers. E.g.,
$\tau = \{1, 1, 3\}$ is a trait in which the datum at index 1 has multiplicity 2, and the datum at index 3 has unit multiplicity.
In our running example, this trait might represent the topic about sports; the first document has two sports words, and the
third document has one sports word.
\bnumdefn \label{defn:trait}
A \emph{trait} is a finite, nonempty multiset of $\nats$.
\enumdefn
Let the set of all traits be denoted $\traits$.
A single trait is not sufficient to capture the combinatorial structure underlying
the first $N\in\nats$ data in the sequence: each datum may be
a member of multiple traits (with varying degrees of membership). The traits have no inherent order
just as the topics ``sports'', ``arts'', and ``science'' have no inherent order.
And each document may contain words from multiple topics.
Building from \cref{defn:trait} and motivated by these desiderata, we define a finite trait allocation
 as a finite multiset of traits.
For example,
$t_4 = \{ \{1\}, \{3, 4\}, \{3, 3\}, \{3, 3\}, \{1, 1, 4\}\}$
represents a collection of traits expressed by the first $4$ data points in a sequence.
In this case, index 1 is a member of two traits, index 2 is a member of none, and so on.
Throughout, we assume that each datum at index $n\in\nats$, $n\leq N$ belongs to only finitely many latent
traits.
Further, for a data set of size $N$,
any index $n>N$ should not belong to any trait; the allocation $t_N$ represents
traits expressed by only the first $N$ data. These statements are formalized in \cref{defn:traitalloc}.
\bnumdefn \label{defn:traitalloc}
A \emph{trait allocation of} $[N]$
is a multiset $t_N$ of traits, where
\[
\forall n \in \nats : n \leq N, \quad &\sum_{\omega\in\traits} t_N(\omega)\cdot \omega(n) < \infty\\
\forall n \in \nats : n > N, \quad &\sum_{\omega\in\traits}t_N(\omega)\cdot \omega(n) = 0.
\]
\enumdefn
Let $\mcT_N$ be the set of trait allocations of $[N]$,
and define $\mcT$ to be the set of all \emph{finite trait allocations}, $\mcT \defined \bigcup_{N}\mcT_N$.
Two notable special cases of finite trait allocations that have appeared in past work
are \emph{feature allocations}~\citep{Griffiths05,Broderick13} and \emph{partitions}~\citep{Kingman78,Pitman95}. Feature allocations are the natural combinatorial
structure underlying feature learning,
where each datum expresses each trait with multiplicity at most 1.
For example, $t_4 = \{ \{1\}, \{3, 4\}, \{3, 1\}, \{3\} \}$ is a feature allocation of $[4]$.
Note that each index may be a member of multiple traits.
Partitions are the natural combinatorial structure underlying clustering,
where the traits form a partition of the indices.
For example,
 $t_4 = \{ \{1, 3, 4\}, \{2\} \}$ is a partition of $[4]$, since its traits
are disjoint and their union is $[4]$.
 The theory in the remainder
of the paper will be applied to recover past results for these structures as corollaries.

Up until this point, we have dealt solely with finite sequences of $N$ data.
However, in many data analysis problems, it is more natural (or at least an acceptable simplifying approximation)
to treat the observed sequence of $N$ data as the beginning of an infinite sequence. As each
datum arrives, it adds its own index to the traits it expresses, and in the process introduces any
previously uninstantiated traits. For example, if after 3 observations we have $t_3 = \{ \{1\}, \{1, 2\}\}$, then
observing the next might yield $t_4 = \{ \{1\}, \{1, 2, 4, 4\}, \{4, 4\}\}$. Note that
when an index is introduced, none of the earlier indices' memberships to traits
are modified; the sequence
of finite trait allocations is \emph{consistent}. To make this rigorous, we define the \emph{restriction} of a trait (allocation), which allows us to relate two trait
allocations $t_N, t_M \in \mcT$ of differing $N$ and $M$.
The restriction operator $\restrict{}{M}$---provided by \cref{defn:restriction} and acting on either traits or finite trait allocations---removes all indices greater than $M$
from all traits, and does not modify the multiplicity of indices
less than or equal to $M$. If any trait becomes empty in this process, it is removed from the
allocation. For example, $\restrict{\{ \{1, 3, 4\}, \{1, 2\}, \{4\} \}}{1} = \{ \{1\}, \{1\} \}$.
Two trait allocations are said to be consistent, per \cref{defn:consistency},
if one can be restricted to recover the other. Thus,
$\{ \{1, 3, 4\}, \{1, 2\}, \{4\} \}$ and $\{ \{1\}, \{1\} \}$ are consistent
finite trait allocations.
\bnumdefn\label{defn:restriction}
The \emph{restriction} $\restrict{}{M} : \traits\to\traits$ of a trait $\tau$ to $M\in\nats$ is defined as
\[
\restrict{\tau}{M}(m) \defined \left\{
                               \begin{array}{ll}
                               \tau(m) & m \leq M\\
                               0 & m > M
                             \end{array}\right. ,
\]
and is overloaded for finite trait allocations $\restrict{}{M} : \mcT \to \mcT_M$ as
\[
\restrict{t_N}{M}(\tau) \defined \left\{
                               \begin{array}{ll}
                              \sum_{\omega \in \traits} \ind(\restrict{\omega}{M} = \tau) \cdot t_N(\omega) & \tau \neq \emptyset\\
                              0 & \tau = \emptyset
                              \end{array}\right. .
\]
\enumdefn
\bnumdefn\label{defn:consistency}
A pair of trait allocations $t_M$ of $[M]$ and $t_N$ of $[N]$ with $M\leq N$ is said to be \emph{consistent}
if $\restrict{t_N}{M} = t_M$.
\enumdefn

The consistency of two finite trait allocations allows us to define the notion of a consistent sequence of trait allocations.
Such a sequence
can be thought of as generated by the sequential process of data arriving; each data point adds its index to its
assigned traits without modifying any previous index.
For example, $\left(\{ \{1\}, \{1\} \},\;\; \{ \{1, 2\}, \{1\} \},\;\; \{ \{1, 2\}, \{1, 3\}\},\;\; \dots\right)$
is a valid beginning to an infinite sequence of trait allocations. The first datum expresses two traits with multiplicity 1,
and the second and third each express a single one of those traits with multiplicity 1.
As a counterexample, $\left(\{ \{1, 1\} \},\;\; \{ \{1, 1\} \},\;\; \{ \{1, 3\} \},\;\; \dots\right)$
is not a valid trait allocation sequence, as the third trait allocation is not consistent with either the first or second.
This sequence does not correspond to building up the traits expressed by data in a sequence; when
the third datum is observed, the traits expressed by the first are modified.
\bnumdefn \label{defn:inftraitseq}
An \emph{infinite trait allocation} $t_\infty = (t_N)$
is a sequence of trait allocations of $[N]$, $N=1, 2, \dots$ for which 
\[\forall N\in\nats, \quad \restrict{t_{N+1}}{N} = t_N.\]
\enumdefn
Note that since restriction is commutative ($\restrict{\restrict{\,\cdot\,}{K}}{M} = \restrict{\restrict{\,\cdot\,}{M}}{K} = \restrict{\,\cdot\,}{K}$ for $K\leq M$),
\cref{defn:inftraitseq} implies that all pairs of elements of the sequence $(t_N)$ are consistent.
Restriction acts on infinite trait allocations in a straightforward way: given $t_\infty = (t_N)$, restriction to $M\in\nats$
is equivalent to the corresponding projection, $\restrict{t_\infty}{M} \defined t_M$.
 
Denote the set of all infinite trait allocations $\mcT_\infty \subset \bigtimes_{N} \mcT_N$. 
Recall that the motivation for developing infinite trait allocations is to capture
the latent combinatorial structure underlying a sequence of observed data. Since this sequence is
random, its underlying structure may also be,
and thus the next task is to develop a corresponding notion of a random infinite trait allocation.
Given a sequence of probability spaces $\left(\mcT_N, 2^{\mcT_N}, \nu_N\right)$ for $N\in\nats$ with consistent measures $(\nu_N)$,
i.e.
\[
\forall N\in\nats, \quad \nu_N(t_N) = \sum_{t_{N+1} \in \mcT_{N+1}}
\ind\left(\restrict{t_{N+1}}{N} = t_N\right)\cdot
 \nu_{N+1}(t_{N+1}),
\]
the Kolmogorov extension theorem \citep[Theorem 5.16]{Kallenberg97} guarantees the existence of a unique random infinite trait allocation $T_\infty$
that satisfies $T_\infty \in \mcT_\infty$ \as and has
finite marginal distributions equal to the $\nu_N$ induced by restriction, i.e.
\[
\forall N\in\nats, \quad \restrict{T_\infty}{N} \dist \nu_N.
\]

The properties of the random infinite trait allocation $T_\infty$ are intimately related to those of the observed sequence of data it represents.
In many applications, the data sequence has the property that its distribution is invariant to finite permutation of its elements;
in some sense, the order in which the data sequence is observed is immaterial. We expect
the random infinite trait allocation $T_\infty$ associated with such an 
\emph{infinite exchangeable sequence}\footnote{For an introduction to exchangeability and related theory, see \citet{Aldous85}.} 
to inherit a similar property.
As a simple illustration of the extension of permutation to infinite trait allocations, suppose we observe the sequence of data $(x_1, x_2, x_3, \dots)$
exhibiting trait allocation sequence $T_1 = \{ \{1, 1\} \}$, $T_2 = \{ \{1, 1, 2\}, \{2\} \}$, $T_3 = \{  \{1, 1, 2\}, \{2\}, \{3, 3\} \}$, and so on.
If we swap $x_1$ and $x_2$ in the data sequence---resulting in the new sequence $(x_2, x_1, x_3, \dots)$---the
traits expressed by $x_2$ become those containing index $1$, the traits for $x_1$ become those containing index $2$, and the rest are unchanged.
Therefore, the permuted infinite trait allocation
is $T'_1 = \{ \{1\}, \{1\} \}$, $T'_2 = \{ \{2, 2, 1\}, \{1\} \}$, $T'_3 = \{  \{2, 2, 1\}, \{1\}, \{3, 3\} \}$, and so on.
Note that $T'_1$ (resp. $T'_2$) is equal to the restriction to $1$ (resp. $2$) of $T_2$ with permuted indices, while
$T'_N$ for $N\geq 3$ is $T_N$ with its indices permuted.
This demonstrates a crucial point---if the permutation affects only indices up to $M\in\nats$ (there is always such an $M$ for finite permutations),
we can arrive at the sequence of trait allocations for the permuted data sequence in two steps. First, we
permute the indices in $T_M$ and then restrict to $1, 2, \dots, M$ to get the first $M$ permuted finite trait allocations. Then
we permute the indices in $T_N$ for each $N>M$.

To make this observation precise, we let $\pi$ be a \emph{finite permutation}
of the natural numbers, i.e.,
\[
\pi : \nats \to \nats, \quad \pi \text{ is a bijection,} \quad \exists M \in \nats : \forall m > M, \, \pi(m) = m,
\]
and overload its notation to operate on traits and (in)finite
trait allocations in \cref{defn:permutation}. Note that if $\pi$ is a finite permutation,
its inverse $\pi^{-1}$ is also a finite permutation with the
same value of $M\in\nats$ for which $m>M$ implies $\pi(m) = m$.
Intuitively, $\pi$ operates on traits and finite trait allocations by permuting their indices. For example, if
$\pi$ has the cycle $(123)$ and fixes all indices greater than 3,
 then $\pi\{1, 1, 2, 4\} = \{2, 2, 3, 4\}$.

\bnumdefn \label{defn:permutation}
Given a finite permutation of the natural numbers $\pi : \nats\to\nats$
that fixes all indices $m > M$,
the permutation of a trait $\tau$ under $\pi$ is defined as
\[
\pi\tau(m) &\defined \tau\left(\pi^{-1}(m)\right), \label{eq:permutetrait}
\]
the permutation of a trait allocation $t_N$ of $[N]$ under $\pi$ is defined as
\[
\pi t_N (\tau) &\defined t_N\left(\pi^{-1}\tau\right),
\]
and the permutation of an infinite trait allocation $t_\infty$ under $\pi$ is defined as
\[
\pi t_\infty &\defined \left(\restrict{ \left(\pi t_{\max(M, N)}\right)}{N}\right)_{N=1}^\infty.\label{eq:permutseq}
\]
\enumdefn
As discussed above, the definition for infinite trait allocations
ensures that the permuted infinite trait allocation is a consistent sequence that corresponds to rearranging
the observed data sequence with the same permutation.
\cref{defn:permutation} provides the necessary framework for studying \emph{infinite exchangeable trait allocations},
 defined as random infinite trait allocations whose distributions are invariant
 to finite permutation.
\bnumdefn \label{defn:exchtraitalloc}
An \emph{infinite exchangeable trait allocation},  $T_\infty$, is a random infinite trait allocation
such that for any finite permutation $\pi : \mathbb{N}\to\mathbb{N}$,
\[
\pi T_\infty \eqd T_\infty.
\]
\enumdefn
Note that if the random infinite trait allocation is
a random infinite partition/feature allocation almost surely, the notion
of exchangeability in \cref{defn:exchtraitalloc} reduces to earlier
notions of exchangeability for random infinite partition/feature allocations~\citep{Kingman78,Aldous85,Broderick13}.
Exchangeability also has an analogous definition for random \emph{finite} trait allocations,
though this is of less interest in the present work.

\begin{figure}
\includegraphics[width=.8\textwidth]{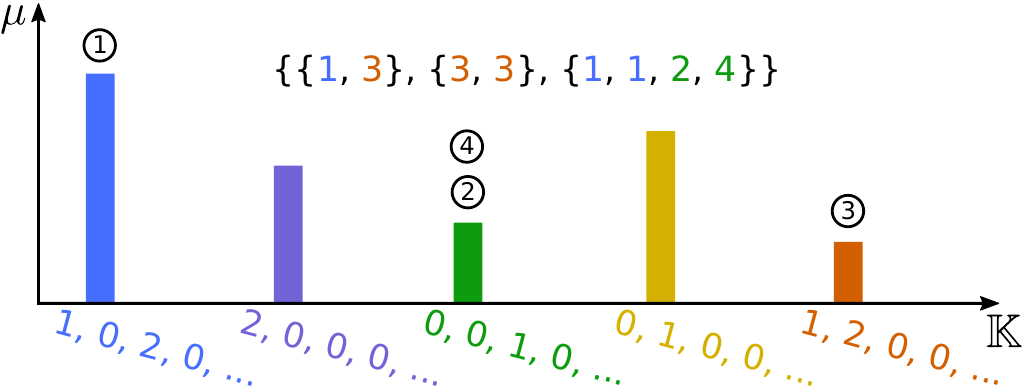}
\caption{An example exchangeable trait allocation construction. For each $N\in\nats$, the trait membership $\xi_N\in\memseqs$
of index $N$ is determined by sampling \iid from the distribution $(\mu_\xi)_{\xi\in\memseqs}$ (depicted by colored bars). The resulting (unordered) trait allocation
for indices up to $4$ is shown above. Here $\xi_1 = (1, 0, 2, 0, \dots)$, $\xi_2=\xi_4=(0,0,1,0,\dots)$, and $\xi_3 = (1,2,0,0, \dots)$.  }\label{fig:egexchtas}
\end{figure}

As a concrete example, consider the
countable set $\memseqs$ of sequences of nonnegative integers $\xi \in (\{0\}\cup\nats)^\infty$ such that $\sum_k \xi_k < \infty$. 
For each data index, we will generate an element of $\memseqs$ and use it to represent a sequence of multiplicities in an ordered sequence of traits.
In particular, we endow $\memseqs$ with probabilities $\mu_{\xi}$ for each $\xi \in \memseqs$.
We start from an empty ordered trait allocation. Then for each data index $N\in \nats$, we sample a sequence $\xi_N \distiid (\mu_\xi)_{\xi\in\memseqs}$; and 
for each $k\in\nats$, we add index $N$ to trait $k$ with multiplicity $\xi_{Nk}$.
The final trait allocation is the unordered collection of nonempty traits.
Since each data index generates its membership in the traits \iid conditioned on $(\mu_\xi)$, the sequence
of trait allocations is exchangeable. This process is depicted in \cref{fig:egexchtas}. 
As we will show in \cref{sec:paintbox}, all infinite exchangeable trait allocations
have a similar construction.

%% file: orderofappearance.tex
\section{Ordered trait allocations and lexicographic ordering}\label{sec:orderofappearance}

We impose no inherent ordering on the traits in a finite trait allocation
via the use of (multi)sets; the allocations $\{\{1\}, \{3, 3\}\}$ and $\{\{3, 3\}, \{1\}\}$
are identical. This correctly
captures our lack of a preferred trait order in many data analysis problems.
However, ordered trait allocations are nonetheless
often useful from standpoints both practical---such as when we need to store a finite trait allocation
in an array in physical memory---and theoretical---such as in developing the characterization of
all  infinite exchangeable trait allocations in \cref{sec:paintbox}.

A primary concern in the development of an ordering scheme is consistency. Intuitively,
as we observe more data in the sequence, we want the sequence of finite ordered trait allocations to ``grow''
but not be ``shuffled''; in other words, if two finite trait allocations are consistent, the traits in their ordered counterparts
at the same index should each be consistent. For partitions, this task is straightforward:
each trait receives as a label its lowest index~\citep{Aldous85}, and the labels are used to order the traits.
This is known as the \emph{order-of-appearance} labeling,
as traits are labeled in the order in which they are instantiated by data in the sequence.
For example, in the partition $t_4 = \{ \{1, 3\}, \{2, 4\}\}$ of $[4]$,
$\{1, 3\}$ would receive label 1 and $\{2, 4\}$ would receive label 2, so $\{1, 3\}$ would be
before $\{2, 4\}$ in the order. Restricting these traits will never change their order---for 
instance, $\restrict{\{1, 3\}}{2} = \{1\}$ and $\restrict{\{2,4\}}{2} = \{2\}$, which
still each receive label 1 and 2, respectively.
If a restriction leaves a trait empty, it is removed and does not interfere with
any traits of a lower label. 
For finite feature allocations, this ordering is inapplicable, since multiple features may have a common lowest index.
Instead, \citet{Griffiths05} introduce a \emph{left-ordered form} in which one feature precedes another 
if it contains an index $n$ that the other does not, and all indices $0< m < n$ have the same membership in both features.
For example, $\{1, 2, 5\}$ precedes $\{1, 3, 5\}$ in this ordering, since the traits both have index $1$, but only 
the first has index $2$.\footnote{Other past work \citep{Broderick13} uses auxiliary randomness to order features, but this technique does not guarantee that 
orderings of two consistent finite trait allocations $t_N, t_M$ are themselves consistent.}
In this section, we show that the well-known \emph{lexicographic ordering}---which generalizes these previous orderings for 
partitions and feature allocations---satisfies our desiderata for an ordering on traits.
 We begin by defining ordered trait allocations.
\bnumdefn \label{defn:orderedtraitalloc}
An \emph{ordered trait allocation} $\ell_N$ of $[N]$ is a sequence $\ell_N = (\ell_{Nk})_{k=1}^K$, $K<\infty$,
of traits $\ell_{Nk} \in \traits$ such that no trait contains an index $n > N$.
\enumdefn
Let $\mcL_N$ be the set of ordered trait allocations of $[N]$, and let $\mcL = \bigcup_{N} \mcL_N$
be the set of all ordered finite trait allocations. 
As in the case of unordered trait allocations, the notion of consistency is intimately tied to that of restriction.
We again require that restriction to $M\in\nats$ removes all indices $m>M$, and removes all traits rendered empty by that process.
However, we also require that the order of the remaining traits is preserved:
for example, if $\ell_3 = (\{3\}, \{1, 2\}, \{2\}, \{1, 1, 2\})$, the restriction of $\ell_3$ to $1$
should yield $(\{1\}, \{1, 1\})$, not $(\{1, 1\}, \{1\})$.
\cref{defn:orderedrestriction} satisfies these desiderata, overloading the $\restrict{}{M}$ function again for notational brevity.
\bnumdefn \label{defn:orderedrestriction}
The \emph{restriction} $\restrict{}{M} : \mcL\to\mcL_M$ of an ordered finite trait allocation $\ell_N$ to $M\in\nats$ is defined as
\[
\restrict{\ell_N}{M} &\defined \squash\left(\left( \restrict{\ell_{Nk}}{M}\right)_{k=1}^K\right),
\]
where the $\squash$ function removes any empty sets from a sequence while preserving the order of the nonempty sets.
\enumdefn
In the example above,
the basic restriction of $\ell_3$ to $1$ would yield $(\emptyset, \{1\},
\emptyset, \{1, 1\})$, which the filter function then
processes to form $\restrict{\ell_3}{1} = \left(\{1\}, \{1, 1\}\right)$, as desired.
Analogously to the unordered case, we say two ordered trait allocations $\ell_N$, $\ell_M$, of $[N]$, $[M]$, with $M\leq N$,
are consistent if $\restrict{\ell_N}{M} = \ell_M$, and define the set of infinite ordered trait allocations $\mcL_\infty$
as the set of infinite sequences of ordered finite trait allocations with $\restrict{\ell_{N+1}}{N}=\ell_N\, \forall N\in\nats$.

Given these definitions, we are now ready to make the earlier intuitive notion
of a consistent trait ordering scheme precise.
\cref{defn:traitallocordering} states that a function $\order{\,\cdot\,}: \mcT\to\mcL$ must
satisfy two conditions to be a valid trait ordering.
The first condition enforces that a trait ordering does not add, remove, or modify the traits in the finite trait allocation $t_N$;
this implies that trait orderings are injective. The second
condition enforces that trait orderings commute with restriction; in other words,
applying a trait ordering to a consistent sequence of finite trait allocations yields a consistent sequence
of ordered finite trait allocations.
For example, suppose
$t_2 =  \{\{2\}, \{1, 2\}\}$,
$t_3 = \{ \{2\}, \{1, 2\}, \{3\}\}$, and
we are given a proposed trait ordering where
$\order{t_2} = \left(\{2\}, \{1, 2\}\right)$ and
$\order{t_3} = \left(\{3\}, \{2\}, \{1, 2\}\right)$.
This would not violate either of the conditions and may be a valid trait ordering.
If instead the ordering was
$\order{t_3} = \left(\{3\}, \{1, 2\}, \{2\}\right)$,
the proposal would not be a valid trait ordering---the
traits $\{2\}$ and $\{1, 2\}$ get ``shuffled'',
i.e., $\order{\restrict{t_3}{2}} = \order{t_2} = \left(\{2\}, \{1, 2\}\right)
\neq
 \left(\{1, 2\}, \{2\}\right) = \restrict{\order{t_3}}{2}$.
\bnumdefn\label{defn:traitallocordering}
A \emph{trait ordering} is a function $\order{\,\cdot\,} : \mcT\to\mcL$ such that:
\begin{enumerate}
\item The ordering is \emph{exhaustive}: If $\order{t_N} = (\tau_k)_{k=1}^K$, then $t_N = \{\tau_1, \dots, \tau_K\}$.
\item The ordering is \emph{consistent}: $\order{\restrict{t_N}{M}} = \restrict{\order{t_N}}{M}$.
\end{enumerate}
\enumdefn

The trait ordering we use throughout is the \emph{lexicographic ordering}: for two traits,
we pick the lowest index with differing multiplicity, and order the one with higher multiplicity first.
For example,
$\{1, 1, 4\} < \{1, 2\}$ since $1$ is the lowest index with differing multiplicity,
and the multiplicity of $1$ is greater in the first trait
than in the second. Similarly,
$\{2, 3\} < \{2, 4\}$ since $3$ has greater multiplicity in the first trait than the second,
and both 1 and 2 have the same multiplicity in both traits.
\cref{defn:traitorder} makes this precise.
\bnumdefn\label{defn:traitorder}
For two traits $\tau, \omega \in \traits$, we say that $\tau < \omega$ if
there exists $n\in\nats$ such that $\tau(n) > \omega(n)$
and all $m\in[n-1]$ satisfy $\omega(m) = \tau(m)$.
\enumdefn
We define $\order{\,\cdot\,} : \mcT \to \mcL$ as the mapping from $t_N$ to
the ordered trait allocation $\ell_N$ induced by the lexicographic ordering. 
The mapping $\order{\,\cdot\,}$ is a trait ordering, as shown by \cref{thm:traitordering}. The proof of \cref{lem:singleconsistent} is provided in \cref{app:proofs}.
\bnlem\label{lem:singleconsistent}
For any pair $\tau, \omega\in\traits$, if $\tau \leq \omega$ then $\restrict{\tau}{M} \leq \restrict{\omega}{M}$ for all $M\in\nats$.
\enlem
\bnthm\label{thm:traitordering}
The mapping $\order{\,\cdot\,}$  is a trait ordering.
\enthm
\bprf
$\order{\,\cdot\,}$ is trivially exhaustive:
since the restriction operation $\restrict{\cdot}{M}$ acts identically to individual traits in both ordered and unordered finite trait allocations,
and empty traits are removed, both $\order{\restrict{t_N}{M}}$ and $\restrict{\order{t_N}}{M}$ have the same multiset of traits (albeit in a potentially different order).
The first trait $\tau$ of $\order{t_N}$ satisfies $\tau \leq \omega$ for any $\omega\in\traits$ such that $t_N(\omega) > 0$,
by definition of $\order{\,\cdot\,}$. By \cref{lem:singleconsistent}, this implies that $\restrict{\tau}{M} \leq \restrict{\omega}{M}$
for all $\omega \in t_N$. Therefore, the first trait
in $\restrict{\order{t_N}}{M}$ is the same as the first trait in $\order{\restrict{t_N}{M}}$.
Applying this logic recursively to $t_N$ with $\tau$ removed, the result follows.
\eprf

%% file: paintbox.tex
\section{De Finetti representation of exchangeable trait allocations}\label{sec:paintbox}

We now derive a de Finetti-style representation theorem 
for infinite exchangeable  trait allocations (\cref{defn:exchtraitalloc})
that extends previous results for partitions and feature allocations~\citep{Kingman78,Broderick13}.
It turns out that \emph{all} infinite exchangeable  trait allocations have essentially 
the same form as in the example construction at the end of \cref{sec:traitalloc},
with some additional nuance.

The high-level proof sketch is as follows.
We first use the lexicographic ordering from \cref{sec:orderofappearance} to associate an \iid sequence of uniform
random labels to the traits in the sequence, in the style of \citet{Aldous85}. We collect the multiset of labels for each index 
into a sequence, called the \emph{label multiset sequence}; the consistency of the ordering from \cref{thm:traitordering} implies that this construction is well-defined.
We show that the label multiset sequence itself is exchangeable in the traditional sequential
sense in \cref{lem:yexch}. And we use 
de Finetti's theorem~\citep[Theorem 9.16]{Kallenberg97} to uncover its construction from conditionally \iid random quantities.
Finally, we relate this construction back to the original set of infinite exchangeable trait allocations to arrive at its representation in \cref{thm:paintbox}.
Throughout the remainder of the paper, $T_\infty \defined (T_N)$ is a random infinite trait allocation and $\phi_\infty \defined (\phi_k) \distiid\distUnif(0, 1)$.

As an example construction of the label multiset sequence, suppose we have $T_4 = \{ \{1, 2, 2\}, \{2, 4\}\}$,
and $\phi_\infty \defined (\phi_k) \distiid \distUnif(0, 1)$. The lexicographic ordering
of $T_4$ is $\order{T_4} = \left( \{1, 2, 2\}, \{2, 4\}\right)$. The first trait in the ordering $\{1, 2, 2\}$ 
receives the first label in the sequence, $\phi_1$, and the second trait $\{2, 4\}$ receives the second label, $\phi_2$.
For each index $n\in[4]$, we now collect the multiset of labels to its assigned traits with the same multiplicity.
Index 1 is a member of only the first trait with multiplicity 1, so its label multiset is $\{\phi_1\}$.
Index 2 is a member of the first trait with multiplicity 2 and the second with multiplicity 1, so its label multiset is  $\{\phi_1, \phi_1, \phi_2\}$.
Similarly, for index 3 it is $\emptyset$, and for index 4 it is $\{\phi_2\}$.
Putting the multisets in order (for index 1, then 2, 3, etc.), the  label multiset sequence is therefore $( \{\phi_1\}, \{\phi_1, \phi_1, \phi_2\}, \emptyset,\{\phi_2\}, \dots)$,
where the ellipsis represents the continuation beyond $T_4$ to $T_5$, $T_6$, and so on.
While the $\phi_k$ may be seen as a mathematical convenience for the proof, an alternative interpretation
is that they correspond to trait-specific parameters in a broader Bayesian hierarchical model. Indeed, our proof would hold
for $\phi_k$ from any nonatomic distribution, not just the uniform. In the document modeling example, 
each $\phi_k$ could correspond to a distribution over English words; $\phi_k$ with high mass on ``basketball'', ``luge'',
and ``curling'' could represent a ``sports'' topic. For this reason, we call the $\phi_k$ \emph{labels}.
Let the set of (possibly empty) finite multisets of $(0, 1)$ be denoted $\labels$.
\bnumdefn\label{defn:labelseq}
The \emph{label multiset sequence} $Y_\infty \defined (Y_N)$ of elements $Y_N\in\labels$
corresponding to $T_\infty$ and $\phi_\infty$
is defined by
\[
  Y_N(\phi) \defined \sum_k \ind\left(\phi=\phi_k\right)\cdot\order{T_N}_k(N). \label{eq:ydef}
\]
\enumdefn
In other words, $Y_N$ is constructed by selecting the $N^\text{th}$ component of $T_\infty$, 
ordering its traits $\tau_1, \dots, \tau_K$, and then adding $\tau_k(N)$ copies of $\phi_k$ to $Y_N$
for each $k \in[K]$. Again, the $\phi_k$ can thus be thought of as labels for the traits, and $Y_N$ 
is the multiset of labels representing the assignment of the $N^\text{th}$ datum to its traits (hence
the name \emph{label multiset sequence}). This construction of $Y_\infty$ ensures that the ``same label applies to the same trait''
as $N$ increases: the \as consistency of the ordering $\order{\,\cdot\,}$ introduced in \cref{sec:orderofappearance} immediately implies that
\[
  \forall N\leq M, \quad Y_N(\phi) \eqas \sum_k \ind\left(\phi=\phi_k\right)\cdot\order{T_M}_k(N).  
\]

\cref{defn:labelseq} implicitly creates a mapping, which we denote $\varphi: \mcT_\infty \times (0, 1)^\infty \to \labels^\infty$.
Since the $\phi_k$ are distinct \as,
we can partially invert $\varphi$ to recover the infinite trait allocation $T_\infty$ corresponding to $Y_\infty$ \as via
\[
T_N(\tau) \eqas \ind\left(\forall n>N, \tau(n)=0\right)\cdot\left|\left\{\phi\in(0, 1) : \forall n\leq N, \tau(n)=Y_n(\phi)\right\}\right|. \label{eq:tvarphidefn}
\]
The first term in the product---the indicator function---ensures that $T_N(\tau)$ is nonzero only for traits $\tau\in\traits$
that do not contain any index $n>N$. The second term counts the number of points $\phi\in(0,1)$ 
for which the multiplicities in $\tau$ match those expressed by the label multiset sequence for $n\leq N$.
Thus, there exists another mapping $\tilde\varphi : \labels^{\infty} \to \mcT_\infty$ 
such that 
\[
\tilde\varphi\left( \varphi\left( T_\infty, \phi_\infty\right) \right) \eqas T_\infty. \label{eq:tvarphidefn2}
\]
The existence of the partial inverse $\tilde\varphi$ is a crucial element in the characterization of all distributions
on infinite exchangeable  trait allocations in \cref{thm:paintbox}. 
In particular, it guarantees that the distributions over random infinite trait allocations are in bijection with
the distributions on label multiset sequences $\labels^\infty$, allowing the characterization of those on $\labels^\infty$
(a much simpler space) instead.
As the primary focus of this work is infinite exchangeable trait allocations, 
we therefore must deduce the particular family of distributions on $\labels^\infty$ that are in bijection
with the infinite exchangeable trait allocations on $\mcT_\infty$.


\cref{lem:yexch} shows that this family is, as one might suspect, 
the exchangeable (in the classical, sequential sense) label multiset sequences.
The main result required for its proof is \cref{lem:pitpiy}, which states that permutation of $T_\infty$ essentially results in 
the same permutation of the components of $Y_\infty$, modulo
reordering the labels in $\phi_\infty$. In other words, permuting the data sequence represented by $T_\infty$ leads to the same permutation of $Y_\infty$.
As an example, consider a setting in which
$T_4 = \{ \{1, 3, 4\}, \{2\}, \{2\}\}$, $\phi_\infty = \left(0.5, 0.4, 0.8, \dots\right)$, and thus
$Y_\infty = \left( \{0.5\}, \{0.4, 0.8\}, \{0.5\}, \{0.5\}, \dots\right)$.
For a finite permutation $\pi$, we define $\pi Y_\infty \defined (Y_{\pi^{-1}(N)})$
and $\pi \phi_\infty \defined (\phi_{\pi^{-1}(k)})$, i.e., permutations 
act on sequences by reordering elements.
If we permute the observed data sequence that $T_4$ represents by $\pi = (12)(34)$, this leads to the permutation of the indices in $T_4$ also by $\pi$, resulting in
$\pi T_4 = \{ \{2, 3, 4\}, \{1\}, \{1\}\}$. If we then reorder $\phi_\infty$ with a different permutation $\pi' = (213)$, 
so $\pi'\phi_\infty = \left(0.4, 0.8, 0.5, \dots\right)$,
then the corresponding label multiset sequence is $Y'_\infty = \left(\{0.4, 0.8\}, \{0.5\}, \{0.5\}, \{0.5\}, \dots\right)$. This $Y'_\infty$
is the reordering of $Y_\infty$ by $\pi$, the same permutation that was used to reorder the observed data; the main result of 
\cref{lem:pitpiy} is that a $\pi'$ \emph{always exists} to reorder $\phi_\infty$
such that this is the case. The proof of \cref{lem:pitpiy} may be found in \cref{app:proofs}.
\bnlem\label{lem:pitpiy}
For each finite permutation $\pi$ and infinite trait allocation $t_\infty$, there exists a finite permutation 
$\pi'$ such that 
\[
  \pi \varphi\left(t_\infty, \phi_\infty\right) \eqas \varphi\left(\pi t_\infty,\pi'\phi_\infty \right).
\]
\enlem

\bnlem\label{lem:yexch}
$T_\infty$ is exchangeable
iff $Y_\infty = \varphi(T_\infty, \phi_\infty)$ is exchangeable.
\enlem
\bprf
Fix a finite permutation $\pi$. Then by \cref{lem:pitpiy} there exists a collection of finite permutations $\pi_{T_\infty}$
that depend on $T_\infty$ such that
\[
\pi Y_\infty &\eqas \varphi(\pi T_\infty, \pi_{T_\infty}\phi_\infty).\label{eq:existspi}
\]
If $Y_\infty$ is exchangeable, then using \cref{eq:existspi} and the definition of $\tilde\varphi$ in \cref{eq:tvarphidefn2},
\[
T_\infty \eqas \tilde\varphi(Y_\infty) \eqd \tilde\varphi(\pi Y_\infty) \eqas \pi T_\infty.
\]
If $T_\infty$ is exchangeable, then again using \cref{eq:existspi} and noting that $\phi_\infty$ is a sequence of \iid random variables and hence also exchangeable,
\[
\pi Y_\infty \eqas \varphi(\pi T_\infty, \pi_{T_\infty}\phi_\infty)
\eqd \varphi(T_\infty, \phi_\infty)
= Y_\infty.
\]
\eprf

We are now ready to characterize all distributions on infinite exchangeable trait allocations in \cref{thm:paintbox} using the \emph{de Finetti representation}
provided by \cref{defn:paintbox}.
At a high level, this is a constructive representation involving three steps. 
Recall that $\memseqs$ is the countable set of sequences of nonnegative integers $(\xi_k)$ such that $\sum_k \xi_k < \infty$. 
First, we generate a (possibly random) distribution over $\memseqs^2$, i.e., a sequence $(\mu_{\xi, \xi'})_{\xi,\xi'\in\memseqs}$ of nonnegative reals such that
\[
\sum_{\xi, \xi' \in \memseqs} \mu_{\xi,\xi'} = 1 \quad \text{and}\quad \forall \xi, \xi'\in\memseqs, \quad \mu_{\xi,\xi'} \geq 0.
\]
Next, for each $N\in\nats$, we sample \iid from this distribution, resulting in two sequences $\xi_N, \xi'_N$.
The sequence $\xi_N$ determines the membership of index $N$ in \emph{regular} traits---which may be joined by other indices---and 
$\xi'_N$ determines its membership in \emph{dust} traits---which are unique to index $N$ and will never be joined by any other index. 
In particular, for each $k\in\nats$, index $N$ joins trait $k$ with multiplicity $\xi_{Nk}$;
and for each $j\in\nats$, index $N$ has $\xi'_{Nj}$ additional unique traits of multiplicity $j$.
For example, in a sequence of documents generated by latent topics, one author may write a single document with
a number of words that are never again used by other authors (e.g.~\emph{Jabberwocky}, by Lewis Carroll);
in the present context, these words would be said to arise from a dust topic. Meanwhile, common collections of words 
expressed by many documents will group together to form regular topics.
Finally, we associate each trait with an \iid $\distUnif(0,1)$ label,
construct the label multiset sequence $Y_\infty$, and use our mapping $\tilde\varphi$ to 
collect these results together to form an infinite trait allocation $T_\infty$.
We say a random infinite trait allocation is \emph{regular} if it has no dust traits with probability 1,
and \emph{irregular} otherwise.

\bnumdefn\label{defn:paintbox}
A random infinite trait allocation $T_\infty$ has a \emph{de Finetti representation}
if there exists a random distribution $(\mu_{\xi,\xi'})$ on $\memseqs^2$ 
such that $T_\infty$ has distribution induced by the following construction: 
\benum
\item generate $(\phi_k), (\phi_{Nj\ell}) \distiid\distUnif(0, 1)$ and $(\xi_N, \xi'_N)\distiid (\mu_{\xi,\xi'})$, 
\item for all $N\in\nats$, define the multisets $R_N, D_N, Y_N$ of $(0, 1)$ via
\[
R_N(\phi) &= \sum_{k,j}\ind(\phi=\phi_k, \xi_{Nk} = j)\cdot j \quad \text{(regular traits)}\label{eq:RNpaintbox}\\
D_N(\phi) &= \sum_{j,\ell}\ind(\phi=\phi_{Nj\ell}, \, \ell \leq \xi'_{Nj})\cdot j \quad \text{(dust traits)} \label{eq:DNpaintbox}\\
Y_N(\phi) &= R_N(\phi)+D_N(\phi),
\]
\item assemble the label multiset sequence $Y_\infty=(Y_N)$ and set $T_\infty = \tilde\varphi(Y_\infty)$.
\eenum
\enumdefn
\cref{thm:paintbox} is the main result of this section, which shows that infinite exchangeable trait allocations---both regular and irregular---are precisely those 
which have a de Finetti representation per \cref{defn:paintbox}. The proof of \cref{thm:paintbox} approaches the problem by characterizing the distribution of the exchangeable label multiset sequence
$Y_\infty$.
%
\bnthm\label{thm:paintbox}
 $T_\infty$ is exchangeable iff it has a de Finetti representation.
\enthm 
\bprf
If $T_\infty$ has a de Finetti representation, then it is exchangeable by the fact
that the $\xi_N, \xi'_N$ are \iid random variables. In the other direction,
if $T_\infty$ is exchangeable,
then there is a random label multiset sequence $Y_\infty = \varphi(T_\infty, \phi_\infty)$ which is exchangeable by \cref{lem:yexch}.
Since we can recover $T_\infty$ from $Y_\infty$ via $T_\infty = \tilde\varphi(Y_\infty)$, it suffices to characterize
 $Y_\infty$ and then reconstruct $T_\infty$.

We split $Y_N$ into its regular $R_N$ and dust $D_N$ components---that represent,
respectively, traits that are expressed by multiple data points 
and those that are expressed only by data point $N$---defined for 
$\phi\in(0, 1)$ by
\[
D_N(\phi) &= \left\{\begin{array}{ll}
                      0 & \exists M\neq N : Y_M(\phi)>0\\
                      Y_N(\phi)  & \text{otherwise}
                    \end{array}\right.\\
R_N(\phi) &= Y_N(\phi) - D_N(\phi).\label{eq:RNdefn}
\]
Choose any ordering $(\phi_k)$ on the countable set $\{\phi \in (0, 1) : \sum_{N}R_N(\phi) > 0\}$.
Next, we extract the multiplicities in $R_N$ and $D_N$ 
via the sequences $\xi_N, \xi'_N \in \memseqs$,
\[
\xi'_{Nj} \defined \left|\left\{\phi\in(0, 1) : D_N(\phi) = j\right\}\right| && \xi_{Nk} \defined R_N(\phi_k).
\] 
Note that we can recover the distribution of $Y_\infty$ 
from that of $(\xi_N, \xi'_N)_{N=1}^\infty$ by generating sequences  
$(\phi'_k), (\phi'_{Nj\ell}) \distiid \distUnif(0, 1)$ and using steps 2 and 3 of \cref{defn:paintbox}.
Therefore it suffices to characterize the distribution of $(\xi_N, \xi'_N)_{N=1}^\infty$.
Note that $(\xi_N, \xi'_N)_{N=1}^\infty$
is a function of $Y_\infty$ such that permuting the elements of $Y_\infty$
corresponds to permuting those of $(\xi_N, \xi'_N)_{N=1}^\infty$ in the same way.
Thus since $Y_\infty$ is exchangeable, so is
$(\xi_N, \xi'_N)_{N=1}^\infty$.
And since $(\xi_N, \xi'_N)_{N=1}^\infty$ is a sequence in a Borel space,
de Finetti's theorem~\citep[Theorem 9.16]{Kallenberg97} states that there exists
a directing random measure $\mu$ such that $(\xi_N, \xi'_N) \distiid \mu$.
Since the set $\memseqs^2$ is countable, we can represent $\mu$
with a probability $\mu_{\xi, \xi'}$
for each tuple $(\xi, \xi')\in\memseqs^2$.
\eprf
The representation in \cref{thm:paintbox} 
generalizes de Finetti representations for both clustering (the Kingman paintbox) and feature allocation (the feature paintbox) \citep{Kingman78,Broderick13},
as shown by \cref{cor:partitionpaintbox,cor:featurepaintbox}.
Further, \cref{cor:featurepaintbox} is the first de Finetti representation for feature allocations
that accounts for the possibility of dust features; previous results were limited
to regular feature allocations~\citep{Broderick13}. 
\cref{thm:paintbox} also makes the distinction between regular and irregular trait allocations
straightforward, as shown by \cref{cor:regtrait}. 
\bncor\label{cor:regtrait}
An exchangeable trait allocation $T_\infty$ is regular iff it has a de Finetti representation where
$\mu_{\xi, \xi'} > 0$ implies $\sum_k \xi'_k = 0$.
\encor
\bncor\label{cor:partitionpaintbox}
A partition $T_\infty$ is exchangeable iff it has a de Finetti representation 
where $\mu_{\xi, \xi'} > 0$ implies either
\bitems
\item $\sum_k \xi_k = 1$ and $\sum_k\xi'_k = 0$, or
\item $\sum_k \xi_k = 0$, $\xi'_1 = 1$, and $\sum_k \xi'_k = 1$.
\eitems
%
\encor
\bncor\label{cor:featurepaintbox}
A feature allocation $T_\infty$ is exchangeable
iff it has a de Finetti representation 
where $\mu_{\xi, \xi'}>0$ implies that
\bitems
\item  $\forall k \in \nats$, $\xi_k \leq 1$,  and $\forall j > 1$, $\xi'_j = 0$.
\eitems
%
%
%
\encor

%% file: freq.tex
\section{Frequency models and probability functions}\label{sec:freq}
The set of infinite exchangeable trait allocations encompasses a very expressive class of random infinite trait allocations:
membership in different regular traits at varying multiplicities can be correlated, 
membership in dust traits can depend on membership in regular traits, etc.
While interesting, this generality makes constructing models with efficient posterior inference procedures difficult. 
A simplifying assumption one can make is that given the directing measure $\mu$, the membership of an index in a particular trait 
is independent of its membership in other traits. 
This assumption is often acceptable in practice, and limits the infinite exchangeable trait allocations 
to a subset---which we refer to as \emph{frequency models}---for which efficient inference is often possible.  
Frequency models, as used in the present context, generalize the notion of a feature frequency model \citep{Broderick13} for feature allocations.

At a high level, this constructive representation consists of three steps.
First, we generate random sequences of nonnegative reals $(\theta_{kj})$ and $(\theta'_j)$
such that $\sum_{k,j}\theta_{kj} < \infty$, 
$\sum_j \theta'_j < \infty$, and $\forall k \in \nats$, $\sum_j \theta_{kj} \leq 1$. 
The quantity $\theta_{kj}$ is the probability that an index joins regular trait $k$ with multiplicity $j$,
while $\theta'_j$ is the average number of dust traits of multiplicity $j$ for each index.
Next, each index $N\in\nats$ independently samples its multiplicity $\xi_{Nk}$ in regular trait $k$ from 
the discrete distribution $(\theta_{kj})_{j=0}^\infty$, where $\theta_{k0}\defined 1-\sum_j \theta_{kj}$
is the probability that the index is not a member of trait $k$. 
For each $j\in\nats$,
each index $N\in\nats$ is a member of an additional $\xi'_{Nj}\distind \distPoiss(\theta'_j)$ dust traits of multiplicity $j$. 
Finally, we collect these results together to form an infinite trait allocation $T_\infty$.
Note that the above essentially imposes a particular form for $\mu$, as given by \cref{defn:freq}.

\bnumdefn\label{defn:freq}
A random infinite trait allocation $T_\infty$ has a \emph{frequency model}
if there exist two random sequences $(\theta_{kj})$, $(\theta'_j)$
of nonnegative real numbers such that
$T_\infty$ has a de Finetti representation with
\[
\mu_{\xi, \xi'} = \left(\prod_{k=1}^\infty \theta_{k\xi_k}\right)\cdot\left( \prod_{j=1}^\infty \frac{(\theta'_j)^{\xi'_j} e^{-\theta'_j}}{\xi'_j!}\right).
\]
%
\enumdefn

Although considerably simpler than general infinite exchangeable trait allocations,
this representation still involves a potentially infinite sequence of parameters;
a finitary representation would be more useful for computational purposes.
In practice, the marginal distribution of $T_N$ provides such a representation~\citep{Griffiths05,Thibaux2007,James14,Broderick15b}.
So rather than considering a simplified class of de Finetti representations, we can alternatively consider a simplified
class of marginal distributions for $T_N$. In previous work on feature allocations \citep{Broderick13}, the analog of
frequency models was shown to correspond to those marginal distributions that depend only on the unordered feature sizes (the so-called
\emph{exchangeable feature probability functions (EFPFs)}). In the following, we develop the generalization of EFPFs 
for trait allocations and show that the same correspondence result holds in this generalized framework.
We let $\kappa(t_N)$ be the number of unique orderings of a trait allocation $t_N$,
\[
\kappa(t_N) \defined \frac{\left(\sum_{\tau\in\traits} t_N(\tau)\right)!}{\prod_{\tau\in\traits}t_N(\tau)!},\label{eq:kappadef}
\]
and use the \emph{multiplicity profile}\footnote{A very similar quantity is known in the population genetics literature
as the \emph{site} (or \emph{allele}) \emph{frequency spectrum}~\citep{Bustamante01}, though it is typically defined there as 
an ordered sequence or vector rather than as a multiset.}
 of $t_N$, given by \cref{defn:multiplicityprofile}, to capture the multiplicities of indices in its traits.
The multiplicity profile of a trait is defined to be the multiset of multiplicities of its elements, while the multiplicity profile
of a finite trait allocation is the multiset of multiplicity profiles of its traits.
As an example, the multiplicity profile of a trait $\{1, 3, 4, 4, 2, 2, 2, 2\}$
is $\{1, 1, 2, 4\}$, since there are two elements of multiplicity 1, one
element of multiplicity 2, and one of multiplicity 4 in the trait.
If we are given the finite trait allocation $\{ \{1, 1, 2\}, \{2\}, \{3\}, \{3, 3, 3, 3, 1\} \}$, then its multiplicity 
profile is $\{ \{1, 2\}, \{1\}, \{1\}, \{1, 4\} \}$. Note that a multiplicity profile is itself
a trait allocation, though not always of the same indices. Here, the trait allocation is of $[3]$, and its multiplicity profile is a trait allocation of $[4]$.
\bnumdefn\label{defn:multiplicityprofile}
The \emph{multiplicity profile} $\mult{\,\cdot\,}:\traits\to\traits$ of a trait $\tau\in\traits$ 
is defined as 
\[
\mult{\tau}(n) \defined \left| \left\{ m \in \nats : \tau(m) = n\right\}\right|, \label{eq:multtau}
\]
and is overloaded for finite trait allocations $\mult{\,\cdot\,}:\mcT\to\mcT$ as
\[
\mult{t_N}(\xi) \defined \sum_{\tau\in\traits}\ind(\mult{\tau} = \xi)\cdot t_N(\tau).
\]
\enumdefn
We also extend \cref{defn:multiplicityprofile} to ordered trait allocations $\ell_N$, where the multiplicity profile is the ordered
multiplicity profiles of its traits, 
i.e.~$\mult{\ell_N}$ is defined such that $\forall k\in\nats, \, \mult{\ell_N}_k \defined \mult{\ell_{Nk}}$.

The precise simplifying assumption on the marginal distribution of $T_N$ that we 
employ in this work is provided in \cref{defn:etpf}, which generalizes 
past work on exchangeable probability functions~\citep{Pitman95,Broderick13}.
\bnumdefn\label{defn:etpf}
A random infinite trait allocation $T_\infty$ \emph{has an exchangeable trait probability function (ETPF)}
if there exists a function $p : \nats\times\mcT\to\reals_+$ such that for all $N\in\nats$,
\[
\Pr\left(T_N = t_N\right)= \kappa(t_N)\cdot p\left(N, \mult{t_N}\right).\label{eq:etpf}
\]
\enumdefn
One of the primary goals of this section is to relate infinite exchangeable trait allocations
with frequency models to those with ETPFs. 
The main result of this section, \cref{thm:etpf}, shows that these two assumptions are actually equivalent:
any random infinite trait allocation $T_\infty$ that has a frequency model (including those with random $(\theta_{kj})$, $(\theta'_j)$ of arbitrary distribution)
has an ETPF, and any random infinite trait allocation with an ETPF
has a frequency model. Therefore, we are able to use the simple construction of frequency models
in practice via their associated ETPFs.
\bnthm\label{thm:etpf}
$T_\infty$ has a frequency model iff it has an ETPF. 
\enthm
The key to the proof of \cref{thm:etpf} is the \emph{uniformly ordered infinite trait allocation},
defined below in \cref{defn:uniformlyorderedtraitalloc}. Recall that $\mcL_\infty$ is the space of consistent, ordered infinite trait allocations
and that $L_\infty$ denotes an ordering of $T_\infty$.
Here, we develop the \emph{uniform ordering} $L_\infty$: intuitively,
for each $N\in\nats$, $L_{N+1}$ is constructed by inserting the new traits in $T_{N+1}$ relative to $T_N$ 
into uniformly random positions among the elements of $L_{N}$.
This guarantees that $L_N$ is marginally a uniform random permutation of $\order{T_N}$
for each $N\in\nats$, and that $L_\infty$ is a consistent sequence, i.e.~$L_\infty\in\mcL_\infty$. 
There are two advantages to analyzing $L_\infty$ rather than 
$T_\infty$ itself. First, the ordering removes 
 the combinatorial difficulties associated with analyzing $T_\infty$. Second, 
the traits
are independent of their ordering,
thereby avoiding the statistical coupling of the ordering based solely on $\order{\,\cdot\,}$.

The definition of the uniform ordering $L_\infty$ in \cref{defn:uniformlyorderedtraitalloc} is based 
on associating traits with the uniformly distributed \iid sequence $\phi_\infty$,
and ordering the traits based on the order of those values.
To do so, we require a definition of the finite permutation $\pi_n$ 
that rearranges the first $n$ elements of $\phi_\infty$
to be in order and leaves the rest unchanged, known as \emph{the $n^\text{th}$ order mapping} $\pi_n$ of $\phi_\infty$. For example,
if $\phi_\infty = \left(0.4, 0.1, 0.3, 0.2, 0.5, \dots\right)$,
then $\pi_3$ is represented in cycle notation as $(321)$,
and $\pi_3\phi_\infty = \left(0.1, 0.3, 0.4, 0.2, 0.5, \dots\right)$.
The precise formulation of this notion is provided by \cref{defn:ordermapping}.
\bnumdefn\label{defn:ordermapping}
The \emph{$n^\text{th}$ order mapping} $\pi_n : \nats\to\nats$ 
of the sequence $\phi_\infty$ is the finite permutation defined by
\[
\pi_n(k) \defined \left\{\begin{array}{ll}
\left|\left\{j\in\nats : j\leq n, \, \phi_j\leq \phi_k\right\}\right| & k \leq  n\\
        k & k > n
\end{array}\right. .
\]
\enumdefn
\cref{defn:uniformlyorderedtraitalloc} shows how to use the $n^\text{th}$ order mapping
to uniformly order an infinite trait allocation:
we rearrange the lexicographic ordering of $T_N$ using the $K_N^\text{th}$ order mapping $\pi_{K_N}$ where
$K_N$ is the number of traits in $T_N$.
\bnumdefn\label{defn:uniformlyorderedtraitalloc}
The \emph{uniform ordering} $L_\infty\defined (L_N)$ of  $T_\infty$
is
\[
\L_{N k} \defined \order{T_N}_{\rho_{N}^{-1}(k)},
\] 
where $\rho_N \defined \pi_{K_N}$ and $K_N = \sum_{\tau\in\traits}T_N(\tau)$ is the number of traits in $T_N$.
\enumdefn
Note that we can also define the \emph{uniformly ordered label multiset sequence} $Y_\infty = (Y_N) \in \labels^\infty$
from the uniform ordering $L_\infty$ of $T_\infty$ via
\[
Y_{N}(\phi) \defined \sum_{k}L_{Nk}(N)\cdot \ind\left(\phi=\phi_{\rho_N^{-1}(k)}\right),
\]
and recover the original infinite random trait allocation $T_\infty \eqas \tilde\varphi(Y_\infty)$ from 
the mapping $\tilde\varphi$ in \cref{eq:tvarphidefn2}.

The proof of \cref{thm:etpf} relies on \cref{lem:Lcondind}, a collection of two technical results associated with 
uniformly ordered infinite trait allocations $L_\infty$ for which the associated unordered infinite trait allocation $T_\infty$ 
has an ETPF. The first result states that 
$L_N$ and $\mult{L_{N+k}}$ are conditionally independent given $\mult{L_N}$ for any $N, k\in\nats$;
essentially, if the distribution of $L_N$ depends only on its multiplicity profile, knowing the multiplicity profiles of further
uniformly ordered trait allocations in the sequence $L_\infty$ does not provide any extra useful information about $L_N$.
The second result states that the distribution of $L_N$ conditioned on $\mult{L_N}$ is uniform.
The proof of \cref{lem:Lcondind} may be found in \cref{app:proofs}. 
\bnlem\label{lem:Lcondind}
If $T_\infty$ has an ETPF, and $L_\infty$ is the uniform ordering of $T_\infty$, then 
for all $N\in\nats$, $\ell_N \in \mcL_N$,
\[
\Pr\left(L_N = \ell_N \given \mult{L_N}, \mult{L_{N+1}}, \mult{L_{N+2}},\dots\right) = \Pr\left(L_N = \ell_N \given \mult{L_N}\right) \quad a.s.,
\]
and $\Pr\left(L_N = \cdot \given \mult{L_N}\right)$ is a uniform distribution over the ordered trait allocations of $[N]$ consistent with $\mult{L_N}$.
\enlem

\bprfof{\cref{thm:etpf}}
Let $L_\infty \defined (L_N)$ be the uniform ordering of $T_\infty \defined (T_N)$.
For any $N\in\nats$, $\ell_N\in\mcL_N$, and $t_N\in\mcT_N$ such that
$\ell_N$ is an ordering of $t_N$, 
\[
\Pr\left(L_N = \ell_N\right) &= \sum_{t'_N \in \mcT_N} \Pr\left(L_N = \ell_N \given T_N = t'_N\right)\Pr\left(T_N = t'_N\right)\\
&=\Pr\left(L_N = \ell_N \given T_N = t_N\right)\Pr\left(T_N = t_N\right)\\
&= \kappa(t_N)^{-1}\Pr\left(T_N = t_N\right),
\]
where the sum collapses to a single term since $t_N\in\mcT_N$ is the unique unordered version of $\ell_N$,
and $\Pr\left(L_N = \ell_N \given T_N = t_N\right) = \kappa(t_N)^{-1}$ since $L_N$ is uniformly distributed over the possible orderings of $T_N$.
 Thus
\[
\Pr\left(T_N = t_N\right) &= \kappa(t_N)\cdot \Pr\left(L_N = \ell_N\right).
\]
Suppose $T_\infty$ has a frequency model as in \cref{defn:freq}. 
To show $T_\infty$ has an ETPF, it remains to show that there exists a function $p$ such that
\[
\Pr\left(L_N = \ell_N\right) &= p\left(N, \mult{t_N}\right).
\]
The major difficulty in doing so is that there is ambiguity in how $L_N=\ell_N$ was
generated from the frequency model;
any trait $\ell_{Nk}$ for which $\mult{\ell_{Nk}}$ is a singleton
(i.e., $\ell_{Nk}$ contains a single unique index)
may correspond to \emph{either} a dust or regular trait. 
Therefore, we must condition on both the frequency model parameters and the (random) dust/regular assignments of the $K$ traits in $\ell_N$.
We let $A_j\subset [K]$, $j\in\nats$ be the set of components of $\ell_N$ corresponding to dust traits of multiplicity $j$.
We further let $Q$ be the set of sequences $(A_j)$ such that
$k\in A_j \implies \mult{\ell_{Nk}} = \{j\}$ for all $k, j \in \nats$, i.e., those that
are possible dust/regular assignments of the traits given $\ell_N$. Note in particular that $Q$ is a function of $\mult{\ell_N}$ but not $\ell_N$.
Then by the tower property,
\[
\Pr(L_N = \ell_N)
&= \EE\left[\Pr(L_N = \ell_N \given (A_j), (\theta_{kj}), (\theta'_j)) \right].\label{eq:expectationLN}
\]
Expanding the inner conditional probability, and defining $A = [K]\setminus \bigcup_j A_j$,
\[
\hspace{-.2cm}\Pr(L_N = \ell_N \given \dots) &= 
\frac{\ind( (A_j)\in Q )}{N^{\sum_j\left|A_j\right|}}
\cdot
\prod_{k=1}^\infty\theta_{k0}^{N}
\cdot\!\!\!
\sum_{\begin{subarray}{c}\sigma : A\to\nats\\\sigma\text{ 1-to-1}\end{subarray}}
\prod_{k\in A}\prod_{j=1}^\infty \left(\frac{\theta_{\sigma(k)j}}{\theta_{\sigma(k)0}}\right)^{\mult{\ell_{Nk}}(j)}\!.\label{eq:conditionedLN}
\]
The first term in the product relates to the dust. Given that we know the positions and multiplicities of dust
in $L_N$, the only remaining randomness is in which index expresses each dust trait; and since $L_N$ has a uniformly
random order, the probability of any index expressing dust at an index is $1/N$. The indicator expresses the
fact that the probability of observing $L_N=\ell_N$ is 0 if it is inconsistent with the dust assignments $(A_j)$.
The second and third terms are the sum over the 
probabilities of all  
ways the $(\theta_{kj})$ could have generated the observed regular traits.

Note that the expression in \cref{eq:conditionedLN}
is a function of only $N$ and $\mult{\ell_N}$, and therefore so is $\Pr\left(L_N=\ell_N\right)$ in \cref{eq:expectationLN}.
But since $L_N$ is a uniformly ordered trait allocation,
$\Pr(L_N=\ell_N)$ is invariant to reordering $\ell_N$, so
it is invariant to reordering $\mult{\ell_N}$; 
and since $\ell_N$ is some ordering of the traits in $t_N$,
$\Pr(L_N=\ell_N)$ is a function of only $\mult{t_N}$ and $N$.
Therefore, there exists some function $p$ such that
\[
\Pr(L_N = \ell_N) = p(N, \mult{t_N}),
\]
and $T_\infty$ has an ETPF as required.

Next, assume $T_\infty$ has an ETPF. 
Consider the finite subsequence $(Y_m)_{m=1}^M$ and $\sigma$-algebra $\mcG_N \defined \sigma\left(\rho_N\phi_\infty, \mult{L_N}\right)$, 
where $M \leq N$, and recall that $\rho_N\phi_\infty$ is the $N^\text{th}$ ordering of $\phi_\infty$, $L_N$ is the uniform ordering of $T_N$,
and $\mult{L_N}$ is its multiplicity profile. 
Note that 
\[
&\Pr\left( (Y_m)_{m=1}^M \given \mcG_N\right)\\
 &=  \sum_{\ell_N \in \mcL_N}\!\!\Pr\left( (Y_m)_{m=1}^M \given \rho_N\phi_\infty, L_N = \ell_N\right)\Pr\left(L_N = \ell_N \given \mult{L_N}\right)\\
 &= \sum_{\ell_N \in \mcL_N}\!\!\Pr\left( (Y_m)_{m=1}^M \given \rho_N\phi_\infty, L_N = \ell_N\right)\Pr\left(L_N = \ell_N \given \left(\mult{L_K}\right)_{K=N}^\infty\right)\\
 &= \Pr\left( (Y_m)_{m=1}^M \given \left(\rho_K\phi_\infty, \mult{L_K}\right)_{K=N}^\infty\right)
\]
almost surely, where the steps follow from the law of total probability, 
the measurability of $\mult{L_N}$ with respect to $\sigma\left(L_N\right)$,
\cref{lem:Lcondind}, and
 the measurability of $\rho_{N+K}\phi_\infty$ with respect to $\sigma(\rho_N\phi_\infty)$ for any $K\in\nats$.
Therefore $\Pr\left( (Y_m)_{m=1}^M \given \mcG_N\right)$ is a reverse martingale in $N$, since $\sigma\left(\rho_K\phi_\infty, \mult{L_K}\right)_{K=N}^\infty$ is 
a reverse filtration; so by the reverse martingale convergence theorem \citep[Theorem 6.23]{Kallenberg97},
there exists a $\sigma$-algebra $\mcG$ such that
\[
 \Pr\left( (Y_m)_{m=1}^M \given \mcG_N\right) \convas \Pr\left( (Y_m)_{m=1}^M\given \mcG\right) \quad N\to\infty.
\]
We now study the properties of the limiting distribution. Denoting $Y_{mk}\defined Y_m(\phi_{\rho_N^{-1}(k)})$ for brevity, note that
the uniform distribution of $L_N$ conditioned on $\mult{L_N}$ implies that
\[
\hspace{-.3cm}\Pr\left( Y_{1k} = j \given (Y_m)_{m=2}^M, \mcG_N\right) = \frac{\mult{L_{Nk}}(j) - \sum_{m=2}^M \ind\left(Y_{mk} = j\right)}{N-M+1}, \,\,\, j\in\nats\cup\{0\}
\]
independently across the trait indices $k\in\nats$. Since $\frac{\sum_{m=2}^M \ind\left(Y_{mk} = j\right)}{N-M+1} \convas 0$ as $N\to\infty$,
we have that $Y_1 \indep (Y_m)_{m=2}^M \given \mcG$. By symmetry, $(Y_m)_{m=1}^M$ are conditionally independent given $\mcG$.
Since this holds for all finite subsequences, the result extends to the infinite sequence: $Y_\infty$ is an \iid sequence conditioned on $\mcG$.
It thus suffices to characterize the limit of $\Pr\left(Y_1 \given \mcG_N\right)$.

Define $\mcD_{Nj}$ to be the set of indices for ``dust-like'' traits of multiplicity $j$, 
 and $\mcR_N$ to be the remaining component indices corresponding to nonempty ``regular-like'' traits, 
\[
\mcD_{Nj} &= \left\{k\in\nats : \mult{L_{Nk}}=\{j\}\right\}, \quad j\in\nats\\
\mcR_N &= \left\{k\in\nats : \mult{L_{Nk}}\neq\emptyset\right\} \setminus \cup_j \mcD_{Nj}.
\]
Simulating from $\Pr\left(Y_1\given\mcG_N\right)$ can be performed in two steps. First, 
independently for every $k\in\mcR_N$, we set $Y_{1k}$
to $j\in\nats$ with probability $\mult{L_{Nk}}(j)/N$, and to 0 with probability $1-\sum_j \mult{L_{Nk}}(j)/N$. 
Then for each $j\in\nats$,
we generate $S_j \dist\distBinom(\left|\mcD_{Nj}\right|, 1/N)$, select a
subset of $\mcD_{Nj}$ of size $S_j$ uniformly at random, and set $Y_{1k}$ for each $k$ in the subset to $j$.
Given the almost-sure convergence of $\Pr\left(Y_1\given\mcG_N\right)$ as $N\to\infty$, the first step implies the existence of a countable sequence $(\phi'_k)$ in $(0, 1)$
(a rearrangement of some subset of the sequence $\phi_\infty$)
and sequences of nonnegative reals $(\theta_{kj})_{j=0}^\infty$ such that
\[
 \theta_{kj} = \lim_{N\to\infty}\frac{\mult{L_{Nk}}(j)}{N}, \quad \theta_{k0}=1-\sum_j\theta_{kj}, \quad \Pr\left(Y_1(\phi'_k) = j \given\mcG\right) = \theta_{kj}
\]
 independently across $k\in\nats$.
Using the law of small numbers~\citep[Theorem 4.6]{Ross11} on the binomial distribution for $S_j$ (with shrinking probabilities $1/N$ as $N\to\infty$), 
and the fact that $\phi_\infty\distiid\distUnif(0, 1)$,
the second step implies that there exists a sequence of positive reals $(\theta'_j)$ 
such that
\[
\theta'_j = \lim_{N\to\infty}\left|\mcD_{Nj}\right|/N,
\]
where $Y_1$ additionally has $\distPoiss(\theta'_j)$
unique elements uniformly distributed on $(0, 1)$ with multiplicity $j$.
Finally, $\sum_j \theta_{kj} \leq 1$ by the above construction, and both $\sum_{k,j}\theta_{kj} < \infty$ 
and $\sum_j \theta'_j < \infty$ almost surely, 
since otherwise the second Borel--Cantelli lemma combined with the \iid nature of 
$Y_\infty$ conditioned on $\mcG$ 
would imply that each $Y_n$ is not a finite multiset, 
which contradicts the assumption that any index is 
a member of only finitely many traits almost surely.
Thus $T_\infty = \tilde\varphi(Y_\infty)$ has a frequency model. 
\eprfof

By setting $\theta_{kj} =\theta'_j = 0$ for all $k,j\in\nats : j > 1$,
\cref{thm:etpf} can be used to recover the correspondence between 
random infinite feature allocations with an exchangeable feature 
probability function (EFPF) and those with a feature frequency model, both defined 
in earlier work by \citet{Broderick13}. In the present context, an EFPF
is an ETPF where $p(N, \mult{t_N}) > 0$ only for $t_N$ that are feature allocations. 
These are exactly the $t_N$ for which
$\mult{t_N}$ only contains traits $\tau$ of the 
form $\{1, 1, 1, \dots, 1\}$, i.e., 
$\mult{t_N}(\tau) > 0$ only if $\forall n>1$, $\tau(n)=0$.
\bncor\label{cor:efpf}
A random infinite feature allocation has a feature frequency model
iff 
it has an EFPF.
\encor
For infinite exchangeable partitions, the result is stronger: \emph{all} exchangeable
infinite partitions have an exchangeable
partition probability function (EPPF)~\citep{Pitman95}, 
defined as a summable symmetric function of the partition sizes times $K!$, where $K$ is the number of partition elements.
\cref{thm:etpf} cannot be directly used to recover this result:
no choice of $(\theta_{kj}), (\theta'_j)$ in \cref{defn:freq} 
or $p(N, \mult{t_N})$ in \cref{defn:etpf} guarantees that the resulting $T_\infty$ 
is a partition. The key issue is that in trait allocations with frequency models, 
the membership of each index in the traits is independent across the traits, while in partitions each index
is a member of exactly one trait. In the EPPF, this manifests itself as an indicator function
that tests whether the traits exhibit a partition structure, where no such test exists in the ETPF
(or EFPF, by extension). 

As trait allocations generalize not only partitions, but other combinatorial structures
with restrictions on index membership as well (cf.~\cref{sec:applications}), it 
is of interest to find a generalization of the correspondence between frequency models and ETPFs
that applies to these constrained structures. We thus require a way of extracting
the memberships of a single index in a trait allocation---referred to as its \emph{membership profile}, as
in \cref{defn:membershipprofile}---so that we can check whether it satisfies
constraints on the combinatorial structure. 
For example, if we have the trait allocation $t_4 = \{ \{1, 1, 2\}, \{1, 2, 3\}, \{1\}\}$, then the membership profile of index $1$ is
$\{1, 1, 2\}$, since index $1$ is a member of two traits with multiplicity 1, and one trait with multiplicity 2.
The membership profile of an index may be empty; for example, here the membership profile of index $4$ in $t_4$ is $\emptyset$.
Finally, and crucially, the membership profile for an index does not change
as more data are observed: for an infinite trait allocation $t_\infty\in\mcT_\infty$,
if $\tau$ is the membership profile of index $n$ in $t_N$ for $n\leq N$,
then for all $M\geq N$, $\tau$ is the membership profile of index $n$ in $t_M$.
\bnumdefn\label{defn:membershipprofile}
The \emph{membership profile of index n} in a finite trait allocation $t_N$
is the multiset $\memb{n}{t_N}$ of $\nats$ defined by
\[
\memb{n}{t_N}(j) &\defined \sum_{\tau\in\traits} \ind\left(\tau(n) = j\right) \cdot t_N(\tau).
\]
\enumdefn
Note that $t_N$ is a partition of $[N]$ if and only if $\forall n\in[N]$ 
$\memb{n}{t_N} = \{1\}$, and $\forall n>N$ $\memb{n}{t_N} = \emptyset$. 
Likewise, $t_N$ is a feature allocation of $[N]$ if and only if 
$\forall n\in[N]$ and $j\in \nats : j>1$, we have $\memb{n}{t_N}(j) = 0$, and $\forall n>N$, $\memb{n}{t_N}=\emptyset$.

\cref{defn:cfreq,defn:cetpf} provide definitions of a frequency model and exchangeable probability function
for combinatorial structures with constraints on the membership
profiles that are analogous to the earlier unconstrained versions in \cref{defn:freq,defn:etpf}.
The intuitive connection to these earlier definitions is made through rejection sampling.
First, we define an acceptable set of membership profiles, known as the 
\emph{constraint set} $\mfC \subset \traits\cup\{\emptyset\}$. Then,
for trait allocations with a \emph{constrained exchangeable trait probability function (CETPF)} in \cref{defn:cetpf},
we generate $T_N$ from the associated unconstrained ETPF and check if all indices $n\in[N]$ have membership
profiles falling in $\mfC$. If this check fails, we repeat the process, and otherwise output $T_N$ 
as a sample from the distribution. Likewise, for trait allocations with a \emph{constrained frequency model},
we generate $Y_n$, $n=1, 2, \dots, N$, progressively checking if all the indices in the associated $T_n$, $n=1, 2, \dots, N$ 
have membership profiles in $\mfC$. If any check fails, we repeat the generation of $Y_n$ for that index $n\in\nats$ until 
it passes. We continue this process until we reach $N\in\nats$ and output $T_N$ as a sample from the distribution. To sample $T_\infty$,
we do the same thing but do not terminate the sequential construction at any finite $N\in\nats$.
Constrained frequency models and CETPFs are the combinatorial analogue of restricted nonparametric processes \citep{Williamson13,DoshiVelez17}.
\bnumdefn\label{defn:cfreq}
A random infinite trait allocation $T_\infty$ has a \emph{constrained frequency model}
with constraint set $\mfC\subset \traits\cup\{\emptyset\}$
if it has a frequency model with step (2) 
from \cref{defn:paintbox} replaced by
\benum
\setcounter{enumi}{1}
\item For $N=1, 2, \dots$,
\crefname{enumi}{step}{steps}
\benum
\item generate $Y_N=R_N+D_N$ as in step (2) of \cref{defn:paintbox}, \label{step:cetpfRNDNYN}
\item let $\mult{Y_N}$ be the multiset of $\nats$ defined by
\[
\mult{Y_N}(n) \defined \left|\left\{ \phi \in (0, 1) : Y_N(\phi) = n\right\}\right|,
\]
\item if $\mult{Y_N}\in\mfC$, continue; otherwise, go to \cref{step:cetpfRNDNYN}. \label{step:cetpfreject}
\eenum
\eenum
\enumdefn
Note that in \cref{defn:cfreq}, $\mult{Y_N}$ is precisely the membership profile of index $N$.
That is to say, if we were to construct
$T_\infty$ from $Y_\infty = \left(Y_1, \dots, Y_N, \emptyset, \emptyset, \dots\right)$, 
then $\mult{Y_N} = \memb{N}{T_N}$. Using $\mult{Y_N}$ instead of this
construction simplifies the definition considerably.
\bnumdefn\label{defn:cetpf}
An infinite trait allocation $T_\infty$ has a \emph{constrained exchangeable trait probability function (CETPF)}
with constraint set $\mfC\subset \traits\cup\{\emptyset\}$ if there exists a function $p : \nats\times\mcT\to\reals_+$ such that 
for all $N\in\nats$,
\[
\sum_{t_N\in\mcT_N}\kappa\left(t_N\right)\cdot p\left(N, \mult{t_N}\right) < \infty
\]
and
\[
\Pr\left(T_N = t_N\right) &= \kappa\left(t_N\right)\cdot p\left(N, \mult{t_N}\right)\cdot\prod_{n=1}^N\ind\left(\memb{n}{t_N}\in\mfC\right).
\]
\enumdefn
The extension of \cref{thm:etpf}---a correspondence between random infinite trait allocations $T_\infty$
with constrained frequency models and CETPFs in \cref{defn:cfreq,defn:cetpf}---that applies to constrained combinatorial structures is given by \cref{thm:cetpf}.
\bnthm\label{thm:cetpf}
$T_\infty$ has a constrained frequency model with constraint set $\mfC$ iff it has a CETPF with constraint set $\mfC$.
\enthm
\bprf
Suppose $T_\infty$ has a constrained frequency model with constraint set $\mfC$. For finite $N\in\nats$,
generating $T_N$ from the constrained frequency model is equivalent to generating it from
the associated unconstrained frequency model (i.e., removing the rejection in {\crefname{enumi}{step}{steps}\cref{step:cetpfreject}} of \cref{defn:cfreq}),
and then rejecting $T_N$ if $\prod_{n=1}^N \ind\left(\memb{n}{T_N}\in\mfC\right) = 0$.
Since generating $T_N$ from an unconstrained frequency model implies it has an ETPF by \cref{thm:etpf}---which inherently satisfies the summability condition in \cref{defn:etpf} because
it is itself a probability distribution---and the final 
rejection step is equivalent to multiplying the distribution of $T_N$ by 
$\prod_{n=1}^N \ind\left(\memb{n}{T_N}\in\mfC\right)$ and renormalizing, $T_\infty$ has a CETPF with constraint set $\mfC$.

Next, suppose $T_\infty$ has a CETPF with constraint set $\mfC$. We can reverse the above logic:
since the associated ETPF is summable, we can generate $T_N$ by simulating from the (normalized) ETPF and rejecting
if $\prod_{n=1}^N \ind\left(\memb{n}{T_N}\in\mfC\right) = 0$. The ETPF has an associated frequency model
by \cref{thm:etpf}. Instead of rejecting $T_N$ after generating all $Y_n$, $n=1, 2, \dots, N$, we can 
reject after each index $n\in\nats$ based on progressively constructing $T_n$, $n=1, 2, \dots, N$.
\eprf

We can, of course, recover \cref{thm:etpf} from \cref{thm:cetpf} by setting $\mfC = \traits\cup\{\emptyset\}$. 
But \cref{thm:cetpf} also allows us to recover earlier results---using a novel proof technique---about the correspondence of infinite exchangeable partitions
and partitions with an EPPF in \cref{cor:eppf}. The proof of \cref{cor:eppf} uses the fact that the EPPF is a constrained EFPF; it is noted
that other connections between classes of probability functions for clustering and feature allocation have been previously established \citep{Roy14}.
\bncor\label{cor:eppf}
An infinite partition $T_\infty$ is exchangeable iff it has an EPPF.
\encor
\bprf
Suppose $T_\infty$ has an EPPF. The EPPF is a CETPF with $\mfC = \{ \{1\} \}$, and thus $T_\infty$ is exchangeable 
by inspection of \cref{defn:cetpf}; the probability
is invariant to finite permutations of the indices.
In the other direction, if $T_\infty$ is an infinite exchangeable partition, then it has a de Finetti representation of the form specified in \cref{cor:partitionpaintbox}; for notational brevity 
define $w_k = \mu_{\xi, \xi'}$ when $\xi_k = 1$ and $w_0 = \mu_{\xi, \xi'}$ when $\xi'_1 = 1$. Note in particular that $\sum_{k=0}^\infty w_k = 1$, and each index $n\in\nats$
selects its trait from the distribution $(w_k)_{k=0}^\infty$, where selecting $0$ implies selecting a dust (or unique) trait.
We seek a constrained frequency model equivalent to this de Finetti representation, so we set 
$\theta_{kj}=\theta'_j=0$ for all $k,j\in\nats : j>1$ and seek $(\theta_{k1})$ and $\theta'_1$  such that
\[
e^{-\theta'_1}\theta'_1 \prod_k \theta_{k0}\propto w_0 
\quad\text{and}\quad
\forall k\in\nats, \, e^{-\theta'_1}\theta_{k1} \prod_{\ell\neq k} \theta_{\ell0}\propto w_k. 
\]
Dividing by $\prod_k\theta_{k0}$, this is equivalent to finding $(\theta_{k1})$ and $\theta'_1$  such that
\[
\theta'_1 \propto w_0 
\quad\text{and}\quad
\forall k\in\nats, \, \frac{\theta_{k1}}{\theta_{k0}}\propto w_k. 
\]
We have a degree of freedom in the proportionality constant, so set 
that equal to 1 and solve each equation by noting that $\theta_{k1}+\theta_{k0}=1$, yielding
\[
\theta_{k1} = \frac{w_k}{w_k+1}\text{ for }k\in\nats, \qquad \theta'_1 = w_0.
\]
The infinite exchangeable partition $T_\infty$ has a constrained frequency model
with constraint set $\mfC = \{\{1\}\}$ based on  $(\theta_{kj})$, $(\theta'_j)$. By \cref{thm:cetpf} it thus has a CETPF 
with the same constraint set $\mfC$, which is an EPPF.
\eprf

%% file: applications.tex
\section{Application: vertex allocations and edge-exchangeable graphs}\label{sec:applications}
%
%
%

A natural assumption for random graph sequences with $\nats$-labeled vertices---arising from online social networks,
protein interaction networks, co-authorship networks,
email communication networks, etc.~\citep{Goldenberg10}---is
that the distribution is projective and invariant to reordering
the vertices, i.e., the graph is \emph{vertex exchangeable}.
 Under this assumption, however,
the Aldous--Hoover theorem~\citep{Aldous81,Hoover79} for exchangeable arrays
guarantees that the resulting graph is either dense or empty almost surely,
an inappropriate consequence when modeling the sparse networks that occur in most
applications~\citep{Mitzenmacher03,Newman05,Clauset09}. Standard statistical models, which are traditionally vertex exchangeable~\citep{Lloyd12},
are therefore misspecified for modeling real-world networks.
This model misspecification has motivated the development and study of a number of 
projective, exchangeable network models that do not preclude
sparsity~\citep{Caron14,Veitch15,Borgs16,Crane16,Cai16,Herlau16,Williamson16}.
One class of such models assumes the network is generated by an exchangeable sequence of (multisets of) edges---the so-called \emph{edge-exchangeable} models~\citep{Broderick15a,Crane15,Cai16,Crane16,Williamson16}.
These models were studied in the generalized hypergraph setting in concurrent work by \citet{Crane16b}.
In this section we provide an alternate view of
edge-exchangeable multigraphs as a subclass of infinite exchangeable  trait allocations
called \emph{vertex allocations}, thus guaranteeing a de Finetti representation.
We also show that the \emph{vertex popularity model}, a standard example of an edge-exchangeable model, is a constrained
frequency model per \cref{defn:cfreq}, thus guaranteeing the existence of a CETPF
which we call the \emph{exchangeable vertex probability function (EVPF)}.
We begin by considering multigraphs without loops, i.e., edges can occur with
multiplicity and all edges contain exactly two vertices.
We then discuss the generalization to
multigraphs with edges that can contain one vertex (i.e., a loop) or two or
(finitely many) more vertices (i.e., a hypergraph).

\begin{figure}
\includegraphics[width=\textwidth]{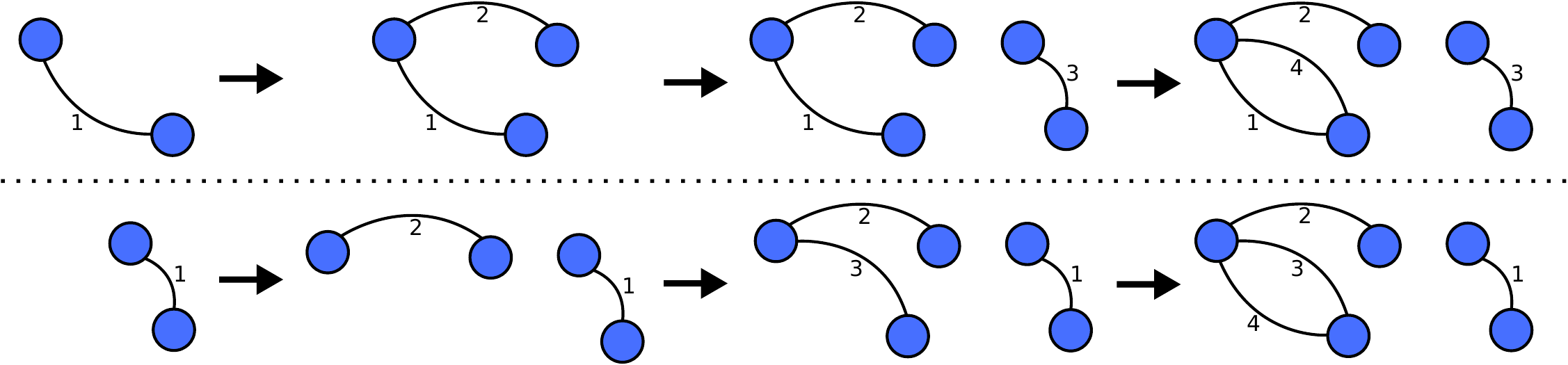}
\caption{Top: the graph encoded by the vertex allocation
$t_4 = \{\{1, 2, 4\}, \{2\}, \{1, 4\}, \{3\}, \{3\} \}$.
The four steps show the sequential construction process of the graph.
Edge labels correspond to indices, and each trait is a vertex.
One or both of the vertices connected to edge 3 and the vertex connected only to edge 2 may be dust; the remaining two are guaranteed
to be regular as they connect to multiple unique edge labels (i.e.~both 1 and 4).
Bottom: the same graph construction with the edges reordered by the permutation $\pi=(314)(2)$, resulting in the vertex
allocation $\pi t_4 = \{ \{4, 2, 3\}, \{2\}, \{4, 3\}, \{1\}, \{1\}\}$. If the vertex
allocation is exchangeable, these sequences have equal probability.}\label{fig:grapheg}
\end{figure}

In the graph setting, the traits correspond to vertices, and the
data indices in each trait correspond to the edges of the graph.
Each data index has multiplicity 1
in exactly two traits---encoding an edge between two separate vertices---as specified in \cref{defn:vertexallocation}.
\cref{fig:grapheg} shows an example encoding of a graph as a vertex allocation.

\bnumdefn\label{defn:vertexallocation}
A \emph{vertex allocation} of $[N]$ is a trait allocation of $[N]$ in which each index
has membership profile equal to $\{1, 1\}$.
\enumdefn
\cref{defn:vertexallocation,thm:paintbox} together immediately yield a de Finetti representation
for edge-exchangeable graphs, provided by \cref{cor:graphpaintbox}.
There are three cases: an edge is either a member of two regular vertices,
one dust vertex and one regular vertex, or
two dust vertices. These three cases are listed in order in \cref{cor:graphpaintbox}.
\bncor\label{cor:graphpaintbox}
An infinite vertex allocation $T_\infty$ is exchangeable iff
it has a de Finetti representation such that
$\mu_{\xi, \xi'} > 0$ implies that either
\benum
\crefname{enumi}{constraint}{constraints}
\item $\exists k\neq j$ such that $\xi_k=\xi_j=1$, $\sum_k \xi_k = 2$, and $\sum_k\xi'_k = 0$,
\item $\exists k $ such that $\xi_k = 1$, $\sum_k \xi_k = 1$, $\xi'_1 = 1$, and $\sum_k \xi'_k = 1$, or
\item $\sum_k\xi_k = 0$, $\xi'_1 = 2$, and $\sum_k \xi'_k = 2$.
\eenum
%
\encor
\cref{defn:vertexallocation,cor:graphpaintbox} can be modified in a number of ways
to better suit the particular application at hand. For example, if loops are
allowed---useful for capturing, for example, authors citing their own earlier work in a
citation network---the membership profile
of each index can be either $\{1, 1\}$ or $\{1\}$.
This allows indices to be a member of a single trait with multiplicity 1, encoding a
loop on a single vertex. If edges between more than two vertices are allowed---that is,
we are concerned with hypergraphs---then we may repurpose the definition of
a feature allocation,
with associated de Finetti representation in \cref{cor:featurepaintbox}, where we view the features
as vertices. 
If $\nats$-valued weights are allowed on the multigraph edges, they
can be encoded using multiplicities greater than 1.
In this case, the index membership profiles must be of the form $\{j, j\}$
for $j\in\nats$, which encodes an edge of weight $j$. Weighted loops may be similarly obtained by
allowing membership profiles of the form $\{j\}$ for $j\in\nats$.
This might be used, for example, to capture an author citing the same work
multiple times in a single document.
Weighted hypergraphs are trait allocations without any restrictions.

\emph{Vertex popularity models}~\citep{Caron14,Cai16,Crane16,Palla16,Herlau16,Williamson16}%
\footnote{These have appeared in previous work as ``graph frequency models''~\citep{Cai16}
or left unnamed, and the weights $w_k$ are occasionally referred to as ``sociability parameters''~\citep{Caron14,Palla16}.}
 are a simple yet powerful class of network models.
There are a number of different versions, but all share the common feature that 
each vertex is associated with a nonnegative weight representing how likely it is to take
part in an edge. Here we adopt a
particular construction based on a sequence of edges: all (potentially infinitely many)
vertices $k\in\nats$ are associated with a
weight $w_k\in(0, 1)$ such that $\sum_k w_k <\infty$,
and we sample an edge between vertex $k$ and $\ell$ with probability proportional
to $w_kw_\ell$.
For an edge-exchangeable vertex popularity model, assuming no loops,
\cref{thm:etpf} enforces that this model has an associated
\emph{exchangeable vertex probability function (EVPF)},
given by \cref{defn:evpf}.
\bnumdefn\label{defn:evpf}
An \emph{exchangeable vertex probability function (EVPF)} is
a CETPF with constraint set $\mfC = \{\{1, 1\}\}$.
\enumdefn
\bncor\label{cor:graphetpf}
A regular infinite exchangeable  vertex allocation has a vertex popularity model
iff it has an EVPF.
\encor
\bprf
We use a similar technique to the proof of \cref{cor:eppf}---we seek
a constrained frequency model
(a sequence $(\theta_{kj})$ and set $\mfC$) that corresponds
to the vertex popularity model with weights $(w_i)$, and then use \cref{thm:cetpf}
to obtain a correspondence with a CETPF (and in particular, an EVPF).
We let $\theta'_j = 0$  for all $j\in\nats$,
let $\theta_{kj} = 0$ for  all $j\in\nats:j>1$,
and seek $(\theta_{k1})$ such that
\[
\forall k,\ell\in\nats:k\neq \ell, \, \theta_{k1}\theta_{\ell1}\prod_{m\neq k,\ell} \theta_{m0}\propto w_kw_j. \label{eq:evpfeq1}
\]
Dividing by $\prod_k \theta_{k0}$, and setting the proportionality constant to 1,
\cref{eq:evpfeq1} is equivalent to
\[
\forall k,\ell\in\nats:k\neq \ell, \, \frac{\theta_{k1}}{\theta_{k0}}\frac{\theta_{\ell1}}{\theta_{\ell0}}=w_kw_\ell.\label{eq:evpfeq2}
\]
\cref{eq:evpfeq2} may be solved, noting that $\forall k\in\nats, \, \theta_{k0}+\theta_{k1}=1$, by
\[
\theta_{k1} = \frac{w_k}{1+w_k} \text{ for } k\in\nats.\label{eq:thetawvpm}
\]
Therefore the vertex popularity model with weights $(w_i)$
is equivalent to a constrained frequency model with $\theta_{k1} = w_k/(1+w_k)$ for $k\in\nats$,
$\theta_{kj} = 0$ for $j>1$, $\theta'_j = 0$ for all $j\in\nats$,
and $\mfC=\{\{1, 1\}\}$ as specified above. \cref{thm:cetpf} guarantees
that the vertex popularity model has a CETPF with constraint set $\mfC$,
and likewise that any CETPF with constraint set $\mfC$ yields a vertex popularity model
by inverting the relation in \cref{eq:thetawvpm}.
\eprf

%% file: conc.tex
\section{Conclusions}

In this work, we formalized the idea of trait allocations---the natural extension 
of well-known combinatorial structures such as partitions and feature allocations to data expressing
latent factors with multiplicity greater than one.
We then developed the framework of exchangeable random infinite trait
allocations, which represent the latent memberships of an exchangeable sequence of data.
The three major contributions in this framework are a de Finetti-style representation theorem
for all exchangeable trait allocations, a correspondence theorem between random trait allocations
with a frequency model and those with an ETPF, 
and finally the introduction and study of the constrained ETPF for capturing
random trait allocations with constrained index memberships.
These contributions apply directly to many other combinatorial structures, such as
edge-exchangeable graphs and topic models.

%% file: proofs.tex
\section{Proofs of results in the main text}\label{app:proofs}
\bprfof{\cref{lem:singleconsistent}}
If $\tau = \omega$, then $\restrict{\tau}{M} = \restrict{\omega}{M}$ (and hence $\restrict{\tau}{M}\leq\restrict{\omega}{M}$) for
all $M \in \nats$ trivially. Otherwise, $\tau < \omega$.
Let $m\in\nats$ be the minimum index in $\tau$ with $\tau(m)>\omega(m)$.
If $M \geq m$, then $\restrict{\tau}{M}(m) > \restrict{\omega}{M}(m)$ and $\restrict{\tau}{M}(j) = \restrict{\omega}{M}(j)$ for $j<m$,
so $\restrict{\tau}{M} < \restrict{\omega}{M}$ by \cref{defn:traitorder}.
If $M < m$, then $\restrict{\tau}{M} = \restrict{\omega}{M}$, since $\tau(n)=\omega(n)$ for any $n < m$.
Therefore, $\restrict{\tau}{M} \leq \restrict{\omega}{M}$.
\eprfof
\bprfof{\cref{lem:pitpiy}}
We prove the result for nonrandom $\phi_\infty$; since the $\pi'$ we develop does 
not depend on $\phi_\infty$, the result holds for all distinct sequences of labels and thus almost surely for \iid uniform $\phi_\infty$ as in the main text as well.

Suppose $\pi$ fixes indices greater than $N$.
Then using \cref{defn:permutation}, $\pi t_N$ is $t_N$ with indices permuted.
Let $K_N = \sum_{\tau\in\traits}t_N(\tau) < \infty$, the number of traits in $t_N$. 
Then let $\pi'$ be the unique finite permutation that maps the index of each trait $\tau$ in $\order{t_N}$ 
to its corresponding trait $\pi\tau$ in $\order{\pi t_N}$, while preserving monotonicity for any traits
of multiplicity greater than 1. 
Mathematically, $\pi'$ 
fixes all $k>K_N$, sets $\pi\left(\order{t_N}_{\pi'^{-1}(k)}\right) = \order{\pi t_N}_k$
for all $k \in [K_N]$, and satisfies $\pi'(k+1) = \pi'(k)+1$ for all $k\in[K_N-1]$ such that $\order{t_N}_k = \order{t_N}_{k+1}$.
Clearly such a permutation exists because $\pi t_N$ contains the same
traits as $t_N$ with indices permuted by \cref{defn:permutation}, and the permutation is unique because any ambiguity (where $t_N$ contains traits with multiplicity greater than 1)
is resolved by the monotonicity requirement.
The monotonicity requirement also implies that $\pi'$ satisfies the desired ordering condition for all $M\geq N$, i.e.
\[
\forall k, M\in \nats : M\geq N,\quad \pi\left(\order{t_M}_{\pi'^{-1}(k)}\right) = \order{\pi t_M}_k,\label{eq:piprimeforallm}
\]
since if an index $M>N$ disambiguates two traits, the fact that $\pi$ fixes all $M>N$ means that 
these two traits have the same relative order in $\order{t_M}$ and $\order{\pi t_M}$.
Set $y'_\infty = \varphi(\pi t_\infty, \pi'\phi_\infty)$.
By \cref{defn:labelseq,eq:piprimeforallm}, we have
\[
\forall M>N, \quad y'_M(\phi_k) &= \order{\pi t_M}_{\pi'(k)}(M) = \pi\left(\order{t_M}_{k}\right)(M),
\]
and since $\pi$ fixes indices greater than $N$ (in particular $\pi(M) = \pi^{-1}(M)= M$),
\[
 \pi\left(\order{t_M}_{k}\right)(M) = \order{t_M}_{k}(M) = y_M(\phi_k) = y_{\pi^{-1}(M)}(\phi_k),
\]
so $y'_\infty = \pi y_\infty$ at all indices greater than $N$. For the remaining indices, we use \cref{defn:labelseq,defn:permutation},
the consistency of the trait ordering in \cref{defn:traitorder}, and the definition of $\pi'$ in sequence:
\[
\forall M\leq N, \quad y'_M(\phi_k) = \order{\restrict{(\pi t_N)}{M}}_{\pi'(k)}(M) = \order{\pi t_N}_{\pi'(k)}(M) = \pi\left(\order{t_N}_{k}\right)(M).
\]
By the definition of permutations of traits in \cref{eq:permutetrait},
\[
 \pi\left(\order{t_N}_{k}\right)(M) = \order{t_N}_{k}(\pi^{-1}(M)).
\]
Finally, again using the consistency of the trait ordering and the fact that $\pi^{-1}(M) \leq N$,
we recover the definition of an element in the original label multiset sequence,
\[
 \order{t_N}_{k}(\pi^{-1}(M)) = \order{t_{\pi^{-1}(M)}}_{k}(\pi^{-1}(M)) = y_{\pi^{-1}(M)}(\phi_k).
\]
Thus $y'_\infty = \pi y_\infty$ at all indices less than or equal to $N$, and the result follows.
\eprfof
\bprfof{\cref{lem:Lcondind}}
For the first statement of the lemma, we need to show that $L_N$ is independent of 
$\mult{L_{N+1}}, \dots, \mult{L_{N+M}}$ given $\mult{L_N}$ for any $M\in\nats$.
We abbreviate $\mult{L_{N+1}}, \dots, \mult{L_{N+M}}$ with $\mult{L_{N\uparrow M}}$,
and abbreviate statements of probabilities by removing unnecessary equalities going forward, e.g.~we replace $\Pr\left(L_N=\ell_N \dots\right)$ with $\Pr\left(L_N\dots\right)$. 
The fact that $\mult{L_N}$ is a function of $L_N$ 
and Bayes' rule yields
\[
\Pr\left(L_N \given \mult{L_N}, \mult{L_{N\uparrow M}}\right) 
&= \frac{\Pr\left(\mult{L_{N\uparrow M}}\given L_N\right)}{\Pr\left(\mult{L_{N\uparrow M}}\given \mult{L_N}\right)}\Pr\left(L_N\given\mult{L_N}\right),
\]
so we require that
\[
\Pr\left(\mult{L_{N\uparrow M}}\given L_N\right) = \Pr\left(\mult{L_{N\uparrow M}}\given \mult{L_N}\right).
\]
Using the law of total probability,
\[
\Pr\left(\mult{L_{N\uparrow M}}\given L_N\right) &= \sum_{L_{N+M}} \Pr\left(\mult{L_{N\uparrow M}}\given L_{N+M}\right) \Pr\left(L_N \given L_{N+M}\right) \frac{\Pr\left(L_{N+M}\right)}{\Pr\left(L_N\right)}.
\]
Since both $L_N$ and $\mult{L_{N\uparrow M}}$ are functions of $L_{N+M}$,
the first two probabilities on the right hand side are actually indicator functions.
Moreover, knowing $L_N$ and $\mult{L_{N\uparrow M}}$ determines $L_{N+M}$ uniquely,
since the differences between $\mult{L_{N+m}}$ and $\mult{L_{N+m+1}}$ for $m=0, 1, \dots, M$ allow one to build up to $L_{N+M}$ sequentially from $L_N$.
Thus there is a unique value $L^\star_{N+M}$ such that
\[
\Pr\left(\mult{L_{N\uparrow M}}\given L_N\right)&= \frac{\Pr\left(L^\star_{N+M}\right)}{\Pr\left(L_N\right)},
\]
where $L^\star_{N+M}$ satisfies
\[
 \restrict{L^\star_{N+M}}{N}=L_N \quad\text{and}\quad \forall m\in[M],\quad \mult{\restrict{L^\star_{N+M}}{N+m}}=\mult{L_{N+m}}.
\]
If we replace $L_N$ with any $L'_N$ such that $\mult{L_N}=\mult{L'_N}$, we have that the corresponding $L'^\star_{N+M}$ satisfies $\mult{L'^\star_{N+M}} = \mult{L^\star_{N+M}}$.
By the ETPF assumption,  the marginal distributions of $L_N$ and $L_{N+M}$ depend on only their multiplicity profiles, so
\[
 \Pr\left(\mult{L_{N\uparrow M}}\given L_N\right) = \frac{\Pr\left(L^\star_{N+M}\right)}{\Pr\left(L_N\right)} &= \frac{\Pr\left(L'^\star_{N+M}\right)}{\Pr\left(L'_N\right)} = \Pr\left(\mult{L_{N\uparrow M}}\given L'_N\right),
\]
and so summing over all such $L'_N$,
\[
 \Pr\left(\mult{L_{N\uparrow M}}\given L_N\right) &=\frac{1}{\left|\left\{L'_N : \mult{L'_N} =\mult{L_N}\right\}\right|} 
\sum_{L'_N : \mult{L'_N}=\mult{L_N}} \Pr\left(\mult{L_{N\uparrow M}}\given L'_N\right).
\]
Therefore $\Pr\left(\mult{L_{N\uparrow M}}\given L_N\right)$ is a function of only $\mult{L_N}$, as desired.
To show that $L_N \given \mult{L_N}$ has the uniform distribution over all ordered trait allocations $L_N$ 
with the given multiplicity profile, 
\[
\Pr\left(L_N \given \mult{L_N}\right) &\propto \Pr\left(\mult{L_N} \given L_N\right)\Pr\left(L_N\right),
\]
which by the ETPF assumption is a constant for any $L_N$ with the given multiplicity profile, and 0 otherwise.
\eprfof